\documentclass[11pt]{article}
\usepackage{times,mathptm}
\usepackage{amsfonts,amscd,amssymb}
\newtheorem{lemma}{Lemma}[section]
\newtheorem{theorem}[lemma]{Theorem}
\newtheorem{proposition}[lemma]{Proposition}
\newtheorem{corollary}[lemma]{Corollary}

\newtheorem{remark}[lemma]{Remark}

\newenvironment{proof}{{\it Proof.}}{\hfill $ \square $ \vskip 4mm}

\newcommand{\thlabel}[1]{\label{th:#1}}
\newcommand{\thref}[1]{Theorem~\ref{th:#1}}
\newcommand{\selabel}[1]{\label{se:#1}}
\newcommand{\seref}[1]{Section~\ref{se:#1}}
\newcommand{\lelabel}[1]{\label{le:#1}}

\newcommand{\prlabel}[1]{\label{pr:#1}}
\newcommand{\prref}[1]{Proposition~\ref{pr:#1}}
\newcommand{\colabel}[1]{\label{co:#1}}
\newcommand{\coref}[1]{Corollary~\ref{co:#1}}

\newcommand{\eqlabel}[1]{\label{eq:#1}}
\newcommand{\eqref}[1]{(\ref{eq:#1})}

\newcommand{\Hom}{\rm{Hom}\,}

\def\lan{\langle}
\def\ran{\rangle}

\def\text#1{\mbox{{\rm #1}}}

\def\ol{\overline}
\def\ul{\underline}
\def\dul#1{\underline{\underline{#1}}}

\def\ot{\otimes}

\def\doublerightleft#1#2{{\lower.2ex\vbox{
\hbox{${\smash{\mathop{\longrightarrow}\limits^{#1}}}$}\vspace*{-4mm}
\hbox{${\smash{\mathop{\longleftarrow}\limits_{#2}}}$}}}}

\input diagrams

\begin{document}
\title{Frobenius and Mashke type Theorems for Doi-Hopf modules
and entwined modules revisited: a unified approach}
\author{Tomasz Brzezi\'nski\thanks{EPSRC Advanced Research Fellow}\\
Department of Mathematics\\
University of Wales Swansea\\ Swansea SA2 8PP, UK\and
S. Caenepeel \\ Faculty of Applied Sciences\\
Free University of Brussels, VUB\\ B-1050 Brussels, Belgium\and
G. Militaru\thanks{Research supported by the bilateral project
``Hopf algebras and co-Galois theory" of the Flemish and
Romanian governments}\\
Faculty of Mathematics\\ University of Bucharest\\
RO-70109 Bucharest 1, Romania\and
Shenglin Zhu\thanks{Research supported by the bilateral project
``New computational, geometric and algebraic methods applied
to quantum groups and differential operators" of the Flemish and
Chinese governments}\\
Institute of Mathematics\\Fudan University\\
Shanghai 200433, China}
\date{}
\maketitle

\begin{abstract}\noindent
We study when induction functors (and their adjoints) between
categories of Doi-Hopf modules and, more generally, entwined
modules are separable, resp. Frobenius. We present a unified
approach, leading to new proofs of old results by the authors,
as well as to some new ones. Also our methods provide a categorical
explanation for the relationship between separability and Frobenius
properties.
\end{abstract}

\section{Introduction}\selabel{0}
Let $H$ be a  Hopf algebra, $A$ an $H$-comodule algebra, and $C$ an
$H$-module coalgebra. Doi \cite{Doi92} and Koppinen \cite{Koppinen95}
independently introduced unifying Hopf modules, nowadays usually
called Doi-Koppinen-Hopf modules, or Doi-Hopf modules. These are
at the same time $A$-modules, and $C$-comodules, with a certain
compatibility relation. Modules, comodules, graded modules,
relative Hopf modules, dimodules and Yetter-Drinfel'd modules are all
special cases of Doi-Hopf modules. Properties of Doi-Hopf modules
(with applications in all the above special cases) have been studied
extensively in the literature. In \cite{CaenepeelMZ97a}, a Maschke
type theorem is given, telling when the functor $F$ forgetting
the $C$-coaction reflects the splitness of an exact sequence, while
in \cite{CaenepeelMZ97b}, it is studied when this  functor is a
Frobenius functor, this means that its right adjoint $\bullet\ot C$
is at the same time a left adjoint.\\
The two problems look very different at first sight, but the results
obtained in \cite{CaenepeelMZ97a} and \cite{CaenepeelMZ97b} indicate
a relationship between them. The main result of \cite{CaenepeelMZ97a}
tells us that we have a Maschke Theorem for the functor $F$ if
$C$ is finitely generated projective and there exists an $A$-bimodule
$C$-colinear map $A\ot C\to C^*\ot A$ satisfying a certain normalizing
condition. In \cite{CaenepeelMZ97b}, we have seen that $F$ is Frobenius
if $C$ is finitely generated and projective and $A\ot C$ and $C^*\ot A$
are isomorphic as $A$-bimodules and $C$-comodules. This isomorphism
can be described using a so-called $H$-integral, this is an element in
$A\ot C$ satisying a certain centralizing condition. The same $H$-integrals
appear also when one studies Maschke Theorems for $G$, the right adjoint
of $F$ (see \cite{CaenepeelIMZ98}). This connection was not well
understood at the time when \cite{CaenepeelMZ97a} and \cite{CaenepeelMZ97b}
were written. The aim of this paper is to give a satisfactory
explanation; in fact we will present a unified approach to both problems,
and solve them at the same time. We will then apply the same technique
for proving new Frobenius type properties: we will study when the
other forgetful functor forgetting the $A$-action is Frobenius, and when
a smash product $A\#_RB$ is a Frobenius extension of $A$ and $B$.
Let us first give a brief overview of new results obtained after
\cite{CaenepeelMZ97a} and \cite{CaenepeelMZ97b}.\\
1) In \cite{CaenepeelIMZ98} and \cite{CaenepeelIMZ99} the notion of
{\it separable functor} (see \cite{NvBvO}) is used to reprove (and
generalize) the Maschke Theorem of \cite{CaenepeelMZ97a}. In fact
separable functors are functors for which a ``functorial" type of
Maschke Theorem holds. A key result due to Rafael \cite{Raf} and
del R\'{\i}o \cite{R} tells us when a functor having a left (resp. right)
adjoint is separable: the unit (resp. the counit) of the adjunction
needs a splitting (resp. a cosplitting).\\
2) {\it Entwined modules} introduced  in
\cite{Brzezinski1} in the context of noncommutative geometry
generalize Doi-Hopf modules.
The most interesting examples of entwined modules turn out to be special
cases of Doi-Hopf modules,
but, on the other hand, the formalism for entwined modules is
more transparent than the one for Doi-Hopf modules. Many results for
Doi-Hopf modules can be generalized to entwined modules, see e.g.
\cite{Brzezinski2}, where the results of \cite{CaenepeelMZ97a}
and \cite{CaenepeelMZ97b} are generalized to the entwined case.\\
3) In \cite{CaenepeelIM}, we look at separable and Frobenius
algebras from the point of view of nonlinear equations; also here we
have a connection between the two notions: both separable and Frobenius
algebras can be described using normalized solutions of the so-called
FS-equation. But the normalizing condition is different in the
two cases.\\
Let $F:\ {\cal C}\to {\cal D}$ be a covariant functor having a right
adjoint $G$. From Rafael's Theorem, it follows that the separability
of $F$ and $G$ is determined by the natural transformations in
$V=\dul{\rm Nat}(GF,1_{\cal C})$ and $W=\dul{\rm Nat}(1_{\cal D},FG)$.
In the case where $F$ is the functor forgetting the coaction, $V$ and
$W$ are computed in \cite{CaenepeelIMZ98}. In fact $V$ and $W$ can also
be used to decide when $G$ is a left adjoint of $F$. This is what we
will do in \seref{2}; we will find new characterizations for $(F,G)$
to be a Frobenius pair, and we will recover the results in
\cite{CaenepeelMZ97b} and \cite{Brzezinski2}. In \seref{4}, we will
apply the same technique to decide when the other forgetful functor
is Frobenius, and in \seref{5}, we will study when the smash product of
two algebras $A$ and $B$ is a Frobenius extension of $A$ and $B$.
This results in necessary and sufficient conditions for the Drinfel'd double
of a finite dimensional Hopf algebra $H$
(which is a special case of the smash product (see \cite{CaenepeelMZ97c})
to be a Frobenius or separable over $H$.\\
We begin with a short section about separable functors and Frobenius
pair of functors. We will explain our approach in the most classical
situation: we consider a ring extension $R\to S$, and consider the
restriction of scalars functor. We derive the (classical) conditions
for an extension to be separable (i.e. the restriction of scalars functor
is separable), split (i.e. the induction functor is separable), and
Frobenius (i.e. restriction of scalars and induction functors
form a Frobenius
pair). We present the results in such a way that they can be extended
to more general situations in the subsequent Sections.\\
Let us remark at this point that the relationship between Frobenius
extensions and separable extensions is an old problem in the literature.
A classical result, due to Eilenberg and Nakayama, tells us that,
over a field $k$, any separable algebra is Frobenius. Several generalizations
of this property exist; conversely, one can give necessary and sufficient
conditions for a Frobenius extension to be separable (see
\cite[Corollary 4.1]{Kadison1}). For more results and a history of this
problem, we refer to \cite{Beidar1}, \cite{Kadison2} and \cite{Kadison3}.\\
We use the formalism
of entwined modules, as this turns out to be more elegant and more
general than that of the  Doi-Hopf modules; several left-right conventions
are possible and there exists a dictionary between them. In
\cite{CaenepeelMZ97a} and \cite{CaenepeelMZ97b}, we have worked with
right-left Doi-Hopf modules; here we will work in the right-right case,
mainly because the formulae then look more natural.\\
Throughout this paper, $k$ is a commutative ring.
We use the Sweedler-Heyneman notation for
comultiplications and coactions. For the comultiplication $\Delta$ on
a coalgebra $C$, we write
$$\Delta(c)=c_{(1)}\ot c_{(2)}.$$
For a right $C$-coaction $\rho^r$ and a left $C$-coaction $\rho^l$ on
a $k$-module $N$, we write
$$\rho^r(n)=n_{[0]}\ot n_{[1]},~~~\rho^l(n)=n_{[-1]}\ot n_{[0]}.$$
We omit the summation symbol $\sum$.

\section{Separable functors and Frobenius pairs of functors}\selabel{1a}
Let $F:\ {\cal C}\to {\cal D}$ be a covariant functor. Recall
\cite{NvBvO} that $F$ is called a {\it separable functor} if the
natural transformation
$${\cal F}:\ \Hom_{\cal C}(\bullet,\bullet)\to \Hom_{\cal
D}(F(\bullet),(\bullet)),$$
induced by $F$ splits. From \cite{Raf} and \cite{R}, we recall the following
characterisation in the case $F$ has an adjoint.

\begin{proposition}\prlabel{1a.1}
Let $G:\ {\cal D}\to {\cal C}$ be a right adjoint of $F$. Let
$\eta:\ 1_{\cal C}\to GF$ and $\varepsilon:\ FG\to 1_{\cal D}$ be the
unit and counit of the adjunction. Then\\
1) $F$ is separable if and only if there exists
$\nu\in V=\dul{\rm Nat}(GF,1_{\cal C})$ such that $\nu\circ\eta=1_{\cal C}$,
the identity natural transformation on ${\cal C}$.\\
2) $G$ is separable if and only if there exists
$\zeta\in W=\dul{\rm Nat}(1_{\cal D},FG)$ such that
$\varepsilon\circ\zeta=1_{\cal D}$,
the identity natural transformation on ${\cal C}$.
\end{proposition}

The separability of $F$ implies a Maschke type Theorem for $F$:
if a morphism $f\in {\cal C}$ is such that $F(f)$ has a one-sided inverse
in ${\cal D}$, then $f$ has a one-sided inverse in ${\cal C}$.\\
A pair of adjoint functors $(F,G)$ is called a {\it Frobenius pair} if $G$ is
not only a right adjoint, but also a left adjoint of $F$.
The following result can be found in any
book on category theory:
$G$ is a left adjoint of $F$ if and only if there exist
natural transformations $\nu\in V=\dul{\rm Nat}(GF,1_{\cal C})$
and $\zeta\in W=\dul{\rm Nat}(1_{\cal D},FG)$ such that
\begin{eqnarray}
F(\nu_M)\circ \zeta_{F(M)}&=& I_{F(M)},\eqlabel{1a.1.1}\\
\nu_{G(N)}\circ G(\zeta_N)&=& I_{G(N)},\eqlabel{1a.1.2}
\end{eqnarray}
for all $M\in {\cal C}$, $N\in {\cal D}$. In order to decide whether $F$ or
$G$ is separable, or whether $(F,G)$ is a Frobenius pair, one has to
investigate
the natural transformations $V=\dul{\rm Nat}(GF,1_{\cal C})$
and $W=\dul{\rm Nat}(1_{\cal D},FG)$. It often happens that the
natural transformations in $V$ and $W$ are determined by single maps.
In this Section we illustrate this in a classical situation and
recover well-known results. In the coming Sections
more general situations are considered.\\
Let $i:\ R\to S$ be a ring homomorphism, and let
$F=\bullet\ot_R S:\ {\cal M}_R\to {\cal M}_S$
be the induction functor. The restriction of scalars functor
$G:\ {\cal M}_S\to {\cal M}_R$ is a right adjoint of $F$. The unit and
counit of the adjunction are
$$\forall M\in {\cal M}_R, \qquad \eta_M:\ M\to M\ot_R S,
~~~\eta_M(m)=m\ot 1,$$
$$\forall N\in {\cal M}_S,\qquad \varepsilon_N:\ N\ot_R S\to N,
~~~\varepsilon_N(n\ot s)=ns.$$
Let us describe $V$ and $W$. Given $\nu:\ GF\to 1_{{\cal M}_R}$ in $V$, it
is not hard to prove that $\ol{\nu}=\nu_R:\ S\to R$ is left and right
$R$-linear.
Conversely, given an $R$-bimodule map $\ol{\nu}:\ S\to R$,
a natural transformation $\nu\in V$ can be constructed by
$$\forall M\in {\cal M}_R, \qquad \nu_M(m\ot s)=m\ol{\nu}(s).$$
Thus we have
\begin{equation}\eqlabel{1a.1.3}
V\cong V_1=\Hom_{R,R}(S,R).
\end{equation}
Now let $\zeta:\ 1_{{\cal M}_S}\to FG$ be  in $W$.
Then $e=\sum e^1\ot e^2=\zeta_S(1)\in S\ot_R S$ satisfies
\begin{equation}\eqlabel{1a.1.4}
\sum se^1\ot e^2=\sum e^1\ot e^2s,
\end{equation}
for all $s\in S$. Conversely if $e$ satisfies \eqref{1a.1.4}, then we can
recover $\zeta$
$$\forall N\in {\cal M}_S,\qquad \zeta_N:\ N\to N\ot_RS,
~~~\zeta_N(n)=ne^1\ot e^2.$$
In the sequel, we  omit the summation symbol, and write
$e=e^1\ot e^2$, where it is understood implicitely that we have a summation.
So we have
$$W\cong W_1= \{e=e^1\ot e^2\in S\ot_R S~|~se^1\ot e^2=e^1\ot e^2s,~
{\rm for~all}~s\in S\}.$$
Combining all these data, we obtain the following result
(cf. \cite{NvBvO} for 1) and 2) and \cite{CaenepeelIM} for 3))

\begin{theorem}\thlabel{1a.2}
Let $i:\ R\to S$ be a ringhomomorphism, $F$ the induction functor, and $G$
the restriction of scalars functor.\\
1) $F$ is separable if and only if there exists a conditional
expectation, that is $\ol{\nu}\in V_1$ such
that $\ol{\nu}(1)=1$, i.e. $S/R$ is a split extension.\\
2) $G$ is separable if and only if there exists a separability idempotent,
that is $e\in W_1$ such that
$e^1e^2=1$, i.e. $S/R$ is a separable extension.\\
3) $(F,G)$ is a Frobenius pair if and only if there exist $\ol{\nu}\in V_1$
and $e\in W_1$ such that
\begin{equation}\eqlabel{1a.2.1}
\ol{\nu}(e^1)e^2=e^1\ol{\nu}(e^2)=1.
\end{equation}
\end{theorem}

\thref{1a.2}~2) explains the terminology for separable functors.
\thref{1a.2}~3) implies the following

\begin{corollary}\colabel{1a.3}
We use the same notation as in  \thref{1a.2}. If $(F,G)$ is a Frobenius pair,
then $S$ is finitely generated and projective as a (right) $R$-module.
\end{corollary}

\begin{proof}
For all $s\in S$, we have
$s=se^1\ol{\nu}(e^2)=e^1\ol{\nu}(e^2s)$, hence $\{e^1,\ol{\nu}(e^2\bullet)\}$
is a dual basis for $S$ as a right $R$-module.
\end{proof}

We have a similar property if $G$ is separable. For the proof
we refer to \cite{Pierce}.

\begin{proposition}\prlabel{1a.4}
With the same notation as in  \thref{1a.2}, if $S$ is an algebra over
a commutative ring $R$,
$S$ is projective as an $R$-module and $G$ is separable, then
$S$ is finitely generated as an $R$-module.
\end{proposition}

Using other descriptions of $V$ and $W$, we find other criteria for $F$
and $G$ to be separable or for $(F,G)$ to be a Frobenius pair. Let
$\Hom_R(S,R)$ be the set of right $R$-module homomorphisms from $S$ to $R$.
$\Hom_R(S,R)$ is an $(R,S)$-bimodule:
\begin{equation}\eqlabel{1a.4.1}
(rfs)(t)=rf(ts),
\end{equation}
for all $f\in \Hom_R(S,R)$, $r\in R$ and $s,t\in S$.

\begin{proposition}\prlabel{1a.5}
Let $i:\ R\to S$ be a ringhomomorphism and use the notation introduced
above. Then
$$V=\dul{\rm Nat}(GF,1_{\cal C})\cong V_2=\Hom_{R,S}(S, \Hom_R(S,R)).$$
\end{proposition}

\begin{proof}
Define $\alpha_1:\ V_1\to V_2$ as follows: for $\ol{\nu}\in V_1$, let
$\alpha_1(\nu)=\ol{\phi}:\ S\to \Hom_R(S,R))$ be given by
$$\ol{\phi}(s)(t)=\ol{\nu}(st).$$
Given $\ol{\phi}\in V_2$, put
$$\alpha^{-1}(\ol{\phi})=\ol{\phi}(1).$$
We invite the reader to verify that $\alpha_1$ and $\alpha_1^{-1}$ are
well-defined and that they are inverses of  each other.
\end{proof}

\begin{proposition}\prlabel{1a.6}
Let $i:\ R\to S$ be a ringhomomorphism and assume that $S$ is finitely
generated and projective as a right $R$-module. Then, with
the notation introduced above,
$$W=\dul{\rm Nat}(1_{\cal D}, FG)\cong W_2=\Hom_{R,S}(\Hom_R(S,R),S).$$
\end{proposition}

\begin{proof}
Let $\{s_i,\sigma_i~|~i=1,\cdots,m\}$ be a finite dual basis of $S$
as a right $R$-module. Then for all $s\in S$ and $f\in \Hom_R(S,R)$,
$$s=\sum_i s_i\sigma_i(s)~~~{\rm and}~~~f=\sum_i f(s_i)\sigma_i.$$
Define $\beta_1:\ W_1\to W_2$ by $\beta_1(e)=\phi$, with
$$\phi(f)=f(e^1)e^2,$$
for all $f\in \Hom_R(S,R)$. To show that $\phi$ is a left $R$-linear
and right $S$-linear map, take any $r\in R$, $s\in S$ and compute
$$\varphi(fs)=f(se^1)e^2=f(e^1)e^2s=\varphi(f)s,$$
$$\varphi(rf)=\sum rf(e^1)e^2=r\varphi(f).$$
Conversely, for $\varphi\in W_2$ define
$$\beta_1^{-1}(\varphi)=e=\sum_i s_i\ot \varphi(\sigma_i).$$
Then for all $s\in S$
\begin{eqnarray*}
\sum_i s_i\ot \varphi(\sigma_i)s&=&
\sum_i s_i\ot \varphi(\sigma_is)\\
&=& \sum_{i,j}\sum_i s_i\ot \varphi(\sigma_i(ss_j)\sigma_j)\\
&=&\sum_{i,j}\sum_i s_i\sigma_i(ss_j)\ot \varphi(\sigma_j)\\
&=&\sum_j ss_j\ot \sigma_j,
\end{eqnarray*}
i.e., $e\in W_1$. Finally, $\beta_1$ and $\beta_1^{-1}$ are  inverses of
each other since
$$\beta_1(\beta_1^{-1}(\varphi))(f)=
\beta_1(\sum_i s_i\ot \varphi(\sigma_i))(f)=
\sum_i f(s_i)\varphi(\sigma_i)=
\sum_i \varphi(f(s_i)\sigma_i)=\varphi(f),$$
$$\beta_1^{-1}(\beta_1(e))=\sum_i s_i\ot \beta_1(e)(\sigma_i)=
\sum_i s_i\ot \sigma_i(e^1)e^2=\sum_i s_i\sigma_i(e^1)\ot e^2=e.$$
\end{proof}

\begin{theorem}\thlabel{1a.7}
Let $i:\ R\to S$ be a ringhomomorphism. We use the notation introduced
above.\\
1) $F:\ {\cal M}_R\to {\cal M}_S$ is separable if and only if there
exists $\ol{\phi}\in V_2$ such that $\ol{\phi}(1)(1)=1$.\\
2) Assume that $S$ is projective as a right $R$-module.
Then $G$ is separable if and only if $S$ is finitely generated as a right
$R$-module and there exists $\phi\in W_2$ such that
$$\sum_i s_i\phi(\sigma_i)=1.$$
3) $(F,G)$ is a Frobenius pair if and only if $S$ is finitely generated
and projective as a right $R$-module, and $\Hom_R(S,R)$ and $S$
are isomorphic as $(R,S)$-bimodules, i.e. $S/R$ is Frobenius.
\end{theorem}

\begin{proof}
The result is a translation of \thref{1a.2} in terms of $V_2$ and $W_2$,
using \prref{1a.4} (for 2)) and \coref{1a.3} (for 3)). We
prove one implication of 3). Assume that $(F,G)$ is a Frobenius pair.
>From \coref{1a.3}, we know that $S$ is finitely generated and projective.
Let $\nu\in V_1$ and $e\in W_1$ be as in part 3) of \thref{1a.2}, and
take $\ol{\phi}=\alpha_1(\ol{\nu})\in V_2$, $\phi=\beta_1(e)\in W_2$.
For all $f\in\Hom_R(S,R)$ and $s\in S$, we have
$$(\ol{\phi}\circ{\phi})(f)(s)=\ol{\nu}(\phi(f)s)=
\ol{\nu}(f(e^1)e^2s)=
f(e^1)\ol{\nu}(e^2s)=f(se^1)\ol{\nu}(e^2)=f(se^1\ol{\nu}(e^2))=f(t)$$
and
$$(\phi\circ\ol{\phi})(s)=\ol{\phi}(s)(e_1)e_2=\ol{\nu}(se^1)e^2=
\ol{\nu}(e^1)e^2s=s.$$
\end{proof}

\section{Entwined modules and Doi-Hopf modules}\selabel{1}
Let $k$ be a commutative ring, $A$ a $k$-algebra, $C$ a (flat)
$k$-coalgebra, and $\psi:\ C\ot A\to A\ot C$ a $k$-linear map.
We use the following notation, inspired by the Sweedler-
Heyneman notation:
$$\psi(c\ot a)= a_{\psi}\ot c^{\psi}.$$
If the map $\psi$ occurs more than once in the same expression,
we also use $\Psi$ or $\Psi'$ as summation indices, i.e.,
$$\psi(c\ot a)= a_{\Psi}\ot c^{\Psi}= a_{\Psi'}\ot c^{\Psi'}.$$
$(A,C,\psi)$ is called a (right-right) {\it entwining structure} if
the following conditions are satisfied for all $a\in A$ and $c\in C$,
\begin{eqnarray}
&& (ab)_{\psi}\ot c^{\psi}= a_{\psi}b_{\Psi}\ot
c^{\psi\Psi},\eqlabel{1.1.1}\\
&& \varepsilon_C(c^{\psi})a_{\psi}=
\varepsilon_C(c)a,\eqlabel{1.1.2}\\
&& a_{\psi}\ot \Delta_C(c^{\psi})=
 a_{\psi\Psi}\ot c_{(1)}^{\Psi}\ot c_{(2)}^{\psi},
\eqlabel{1.1.3}\\
&& 1_{\psi}\ot c^{\psi}=1\ot c.\eqlabel{1.1.4}
\end{eqnarray}
A $k$-module $M$ together with a right $A$-action and a right $C$-coaction
satisfying the compatibility relation
\begin{equation}\eqlabel{1.1.5}
\rho^r(ma)=m_{[0]}a_{\psi}\ot m_{[1]}^{\psi}
\end{equation}
is called an {\it entwined module}. The category of entwined modules and
$A$-linear $C$-colinear maps is denoted by ${\cal C}={\cal M}(\psi)_A^C$.
An important class of examples comes from
{\sl Doi-Koppinen-Hopf structures}.
A (right-right) Doi-Koppinen-Hopf structure consists of a triple $(H,A,C)$,
where
$H$ is a $k$-bialgebra, $A$ a right $H$-comodule algebra, and $C$ a right
$H$-module coalgebra. Consider the map $\psi:\ C\ot A\to A\ot C$ given by
$$\psi(c\ot a)=a_{[0]}\ot ca_{[1]}.$$
Then $(A,C,\psi)$ is an entwining structure, and the compatibility relation
\eqref{1.1.5} takes the form
\begin{equation}\eqlabel{1.1.6}
\rho^r(ma)=m_{[0]}a_{[0]}\ot m_{[1]}a_{[1]}.
\end{equation}
A $k$-module with an $A$-action and a $C$-coaction satisfying \eqref{1.1.6}
is called a {\it Doi-Koppinen-Hopf module} or a {\it Doi-Hopf module}.
Doi-Koppinen-Hopf modules were introduced independently by Doi in
\cite{Doi92} and
Koppinen in \cite{Koppinen95}. Properties of these modules were studied
extensively during the last decade, see e.g. \cite{CaenepeelIMZ98},
\cite{CaenepeelIMZ99}, \cite{CaenepeelMZ97a}, \cite{CaenepeelMZ97b},
\cite{CaenepeelMZ97c}, \cite{CaenepeelR}. Another class of entwining
structures is related to coalgebra Galois extensions, see
\cite{BrzezinskiHajac} for details.
Entwining structures were
introduced in \cite{BrzezinskiM}. Many properties
of Doi-Hopf modules can be generalized to entwined modules (see
e.g. \cite{Brzezinski2}, \cite{Brzezinski99}). Although
 the most studied examples of entwined modules
(graded modules, Yetter-Drinfel'd modules, dimodules,
Hopf modules) are special cases of Doi-Hopf modules, their
properties  can be formulated
more elegantly
in the language of entwined modules.\\
The functor $F:\ {\cal C}={\cal M}(\psi)_A^C\to {\cal M}_A$ forgetting
the $C$-coaction has a right adjoint $G=\bullet\ot C$. The structure
on $G(M)=M\ot C$ is given by the formulae
\begin{eqnarray}
\rho^r(m\ot c)&=& m\ot c_{(1)}\ot c_{(2)},\eqlabel{1.1.7}\\
(m\ot c)a&=& ma_{\psi}\ot c^{\psi}.\eqlabel{1.1.8}
\end{eqnarray}
For later use, we list
the unit and counit natural transformations describing the adjunction,
$$\rho:\ 1_{\cal C}\to GF~~{\rm and}~~\varepsilon:\ FG\to 1_{{\cal M}_A},$$
$$\rho_M:\ M\to M\ot C,~~~\rho_M(m)=\sum m_{[0]}\ot m_{[1]},$$
$$\varepsilon_N=I_N\ot \varepsilon_C:\ N\ot C\to N.$$
In particular, $A\ot C\in {\cal M}(\psi)_A^C$. $A\ot C$ is also a left
$A$-module, the left $A$-action is given by $a(b\ot c)=ab\ot c$.
This makes $A\ot C$ into an object of ${}_A{\cal M}(\psi)_A^C$, the
category of entwined modules with an additional left $A$-action that
is right $A$-linear and right $C$-colinear.\\
The other forgetful functor $G':\ {\cal M}(\psi)_A^C\to \cal{M}^C$
has a left adjoint $F'=\bullet\ot A$. The structure on $F'(N)=N\ot A$
is now given by
\begin{eqnarray}
\rho^r(n\ot a)&=& n_{[0]}\ot a_{\psi}\ot n_{[1]}^{\psi},\eqlabel{1.7.9}\\
(n\ot a)b&=& n\ot ab.\eqlabel{1.7.10}
\end{eqnarray}
The unit and counit of the adjunction are
$$\mu:\ F'G'\to 1_{\cal C}~~{\rm and}~~\eta:\ 1_{{\cal M}^C}\to G'F',$$
$$\mu_M:\ M\ot A\to A,~~~\mu_M(m\ot a)=ma,$$
$$\eta_N:\ N\to N\ot A,~~~\eta_N(n)=n\ot 1.$$
In particular $G'(C)=C\ot A\in {\cal M}(\psi)_A^C$. The map
$\psi:\ C\ot A\to A\ot C$ is a morphism in ${\cal M}(\psi)_A^C$.
$C\ot A$ is also a left $C$-comodule, the left $C$-coaction being induced
by the comultiplication on $C$. This coaction is right $A$-linear and
right $C$-colinear, and thus $C\ot A$ is an object of
${}^C{\cal M}(\psi)_A^C$, the category of entwined modules together with
a right $A$-linear right $C$-colinear left $C$-coaction.

\section{The functor forgetting the coaction}\selabel{2}
Let $(A,C,\psi)$ be a right-right entwining structure,
$F:\ {\cal M}(\psi)_A^C\to {\cal M}_A$ the functor forgetting the
coaction, and $G=\bullet\ot C$ its adjoint. In \cite{CaenepeelMZ97b}
 necessary and sufficient conditions for $(F,G)$ to
be a Frobenius pair are given (in the Doi-Hopf case;
the results were generalized
to the entwining case  in \cite{Brzezinski2}), under
the additonal assumption that $C$ is projective as a $k$-module. In this
Section we  give an alternative characterization that also holds
if $C$ is not necessarily projective, and we  find a new proof of
the results in \cite{CaenepeelMZ97b} and \cite{Brzezinski2}. The
method of proof is the same as in \seref{1a}, i.e., based on explicit
descriptions of $V$ and $W$. These descriptions can be found in
\cite{CaenepeelIMZ98}, \cite{CaenepeelIMZ99} and \cite{Brzezinski99}
in various degrees of generality. To keep this paper self-contained,
we give a sketch of proof. We first investigate
$V=\dul{\rm Nat}(GF,1_{\cal C})$.
Let $V_1$ be the $k$-module consisting of all
$k$-linear maps $\theta:\ C\ot C\to A$ such that
\begin{eqnarray}
\theta(c\ot d)a&=& a_{\psi\Psi}\theta(c^{\Psi}\ot d^{\psi}),\eqlabel{2.1.3}\\
\theta(c\ot d_{(1)})\ot d_{(2)}&=&
\theta(c_{(2)}\ot d)_{\psi}\ot c_{(1)}^{\psi}.\eqlabel{2.1.4}
\end{eqnarray}

\begin{proposition}\prlabel{2.2}
The map $\alpha:\ V\to V_1$ given by $\alpha(\nu)=\theta$,
with
\begin{equation}\eqlabel{2.2.1}
\theta(c\ot d)=(I_A\ot \varepsilon_C)(\nu_{A\ot C}(1_A\ot c\ot d)),
\end{equation}
is an isomorphism of $k$-modules. The inverse $\alpha^{-1}(\theta)=\nu$ is
defined as follows: $\nu_M:\ M\ot C\to M$
is given by
\begin{equation}\eqlabel{2.2.2}
\nu_M(m\ot c)= m_{[0]}\theta(m_{[1]}\ot c).
\end{equation}
\end{proposition}

\begin{proof}
Consider $\ul{\nu}=\nu_{A\ot C}$ and $\ol{\nu}=\nu_{C\ot A}$.
Due to the naturality of $\nu$ and  \eqref{1.1.1} there is a commutative
diagram
$$\begin{diagram}
C\ot A\ot C&\rTo^{\ol{\nu}}&C\ot A&\rTo^{\varepsilon_C\ot I_A}&A\\
\dTo^{\psi\ot I_C}&&\dTo^{\psi}&&\dTo{I_A}\\
A\ot C\ot C&\rTo^{\ul{\nu}}&A\ot C&\rTo^{I_A\ot\varepsilon_C}&A
\end{diagram}$$
Write $\ol{\lambda}=(I_A\ot\varepsilon_C)\circ \ol{\nu}$ and
$\ul{\lambda}=(\varepsilon_C\ot I_A)\circ \ul{\nu}$. Then it follows
that
$$\theta(c\ot d)=\ol{\lambda}(c\ot 1\ot d)=\ul{\lambda}(1\ot c\ot d).$$
We have seen before that $A\ot C\in {}_A{\cal M}(\psi)_A^C$. It is easy
to prove that $GF(A\ot C)=A\ot C\ot C\in {}_A{\cal M}(\psi)_A^C$ - the left
$A$-action is induced by the multiplication in $A$ - and $\ul{\nu}$ is
a morphism in ${}_A{\cal M}(\psi)_A^C$. Thus $\ul{\nu}$ and $\ul{\lambda}$
are left and right $A$-linear, and
\begin{eqnarray*}
\theta(c\ot d)a&=& \ul{\lambda}(1\ot c\ot d)a
=\ul{\lambda}(a_{\psi\Psi}\ot c^{\Psi}\ot d^{\psi})\\
&=& a_{\psi\Psi}\ul{\lambda}(1\ot c^{\Psi}\ot d^{\psi})
= a_{\psi\Psi}\theta(c^{\Psi}\ot d^{\psi}),
\end{eqnarray*}
proving \eqref{2.1.3}. To prove \eqref{2.1.4}, we first observe that
$C\ot A,~GF(C\ot A)=C\ot C\ot A\in {}^C{\cal M}(\psi)_A^C$, the left
$C$-coaction is induced by comultiplication in $C$ in the first factor.
Also $\ol{\nu}$ is a morphism in ${}^C{\cal M}(\psi)_A^C$, and we
conclude that $\ol{\nu}$ is left and right $C$-colinear. Take $c,d\in C$,
and put
$$\ol{\nu}(c\ot d\ot 1)=\sum_i c_i\ot a_i.$$
Writing down the condition that $\ol{\nu}$ is left $C$-colinear, and then
applying $\varepsilon_C$ to the second factor, we find that
\begin{equation}\eqlabel{2.2.3}
c_{(1)}\ot \theta(c_{(2)}\ot d)=\sum_i c_i\ot a_i=
\ol{\nu}(c\ot d\ot 1).
\end{equation}
Since $\ol{\nu}$ is also right $C$-colinear,
$$\ol{\nu}(c\ot 1\ot d_{(1)})\ot d_{(2)}=\sum_i c_{i(1)}\ot a_{i\psi}\ot
c_{i(2)}^{\psi}$$
and, applying $\varepsilon_C$ to the second factor, we find
\begin{equation}\eqlabel{2.2.4}
\theta(c\ot d_{(1)})\ot d_{(2)}=\psi(\sum_i c_i\ot a_i),
\end{equation}
and \eqref{2.1.4} follows from \eqref{2.2.3} and \eqref{2.2.4}. This
proves
that there is a well-defined map $\alpha:\ V\to V_1$.\\
To show that the map $\alpha^{-1}$ defined by \eqref{2.2.2}
is well-defined, take $\theta\in V_1$, $M\in {\cal C}$, and let $\nu_M$
be given by \eqref{2.2.2}. It needs to be shown that $\nu_M\in {\cal C}$, i.e.,
$\nu_M$ is right $A$-linear and right $C$-colinear,
and that $\nu$ is a natural transformation. The right $A$-linearity
follows from \eqref{2.1.3},
and the right $C$-colinearity from \eqref{2.1.4}. Given any
morphism  $f:\ M\to N$ in ${\cal C}$, one easily checks
that for all $m\in M$ and $c\in C$
$$\nu_N(f(m)\ot c)= f(m_{[0]})\theta(m_{[1]}\ot c)
= f(m_{[0]}\theta(m_{[1]}\ot c))
= f(\nu_M(m\ot c)),$$
i.e., $\nu$ is natural. The verification that $\alpha$
and $\alpha^{-1}$ are inverses of
each other is left to the reader.
\end{proof}

Now we  give a description of $W=\dul{\rm Nat}(1_{{\cal M}_A},FG)$. Let
$$W_1=\{z\in A\ot C~|~az=za,~{\rm for~all~}a\in A\},$$
i.e., $z=\sum_l a_l \ot c_l\in W_1$ if and only if
\begin{equation}\eqlabel{2.3.1}
\sum_l aa_l\ot c_l=\sum_l a_la_{\psi}\ot c_l^{\psi}.
\end{equation}

\begin{proposition}\prlabel{2.3}
Let $(A,C,\psi)$ be a right-right entwining structure. Then there is
an isomorphism of $k$-modules $\beta:\ W\to W_1$ given by
\begin{equation}\eqlabel{2.3.2}
\beta(\zeta)=\zeta_A(1).
\end{equation}
The inverse of $\beta$ is $\beta^{-1}(\sum_l a_l\ot c_l)=\zeta$, with
$\zeta_N:\ N\to N\ot C$
given by
\begin{equation}\eqlabel{2.3.3}
\zeta_N(n)=\sum_l na_l\ot c_l.
\end{equation}
\end{proposition}

\begin{proof}
We leave the details to the reader; the proof relies on the fact that
$\zeta_A$ is left and right $A$-linear.
\end{proof}

In \cite{CaenepeelIMZ98}, Propositions \ref{pr:2.2} and \ref{pr:2.3} are
used to determine when the functor $F$ and its adjoint $G$ are separable.

\begin{theorem}\thlabel{2.4}
Let $F:\ {\cal M}(\psi)^C_A\to {\cal M}_A$ be the forgetful functor,
and $G=\bullet\ot C$ its adjoint.\\
$F$ is separable if and only if there exists $\theta\in V_1$ such that
$$\theta \circ \Delta_C=\varepsilon_C.$$
$G$ is separable if and only if there exists $z=\sum_l a_l\ot c_l\in W_1$
such that
$$\sum_l \varepsilon_C(c_l)a_l=1.$$
\end{theorem}

\begin{proof}
This follows immediately from Propositions \ref{pr:1a.1},
\ref{pr:2.2} and \ref{pr:2.3}
\end{proof}

Next we show that the fact that $(F,G)$ is a Frobenius pair is also
equivalent to the existence of $\theta\in V_1$ and $z\in W_1$, but now
satisfying different normalizing conditions.

\begin{theorem}\thlabel{2.5}
Let $F:\ {\cal M}(\psi)^C_A\to {\cal M}_A$ be the forgetful functor,
and $G=\bullet\ot C$ its adjoint. Then $(F,G)$ is a Frobenius pair if
and only if there exist $\theta\in V_1$ and $z=\sum_l a_l\ot c_l\in W_1$
such that the following normalizing condition holds, for all $d\in C$:
\begin{eqnarray}
\varepsilon_C(d)1&=&\sum_l a_l\theta(c_l\ot d)\eqlabel{2.5.1a}\\
&=&\sum_l a_{l\psi}\theta(d^{\psi}\ot c_l).\eqlabel{2.5.1b}
\end{eqnarray}
\end{theorem}

\begin{proof}
Suppose that $(F,G)$ is a Frobenius pair. Then there exist
$\nu\in V$ and $\zeta\in W$ such that (\ref{eq:1a.1.1}-\ref{eq:1a.1.2})
hold. Let $\theta=\alpha(\nu)\in V_1$, and $z=\sum_l a_l\ot c_l=
\beta(\zeta)\in W_1$. Then \eqref{1a.1.1} can be rewritten as
\begin{equation}\eqlabel{2.5.2}
\nu_M(\sum_l ma_l\ot c_l)=
\sum m_{[0]}a_{l\psi}\theta(m_{[1]}^{\psi}\ot c_l)=m,
\end{equation}
for all $m\in M\in {\cal M}(\psi)_A^C$. Taking $M=C\ot A$, $m=d\ot 1$,
one obtains \eqref{2.5.1b}.\\
For all $n\in N\in {\cal M}_A$ and $\in C$, one has
\begin{eqnarray*}
\nu_{G(N)}(G(\zeta_N)(n\ot d))&=&
\nu_{G(N)}(\sum_l na_l\ot c_l\ot d)\\
&=& \sum_l (na_l\ot c_{l(1)})\theta(c_{l(2)}\ot d)\\
&=& \sum_l na_l\theta(c_{l(2)}\ot d)_{\psi}\ot c_{l(1)}^{\psi}\\
{\rm \eqref{2.1.4}}~~~~&=& \sum_l na_l\theta(c_{l}\ot d_{(1)})\ot d_{(2)}
\end{eqnarray*}
and \eqref{1a.1.2} can be written as
\begin{equation}\eqlabel{2.5.3}
n\ot d= \sum_l na_l\theta(c_{l}\ot d_{(1)})\ot d_{(2)},
\end{equation}
for all $n\in N\in {\cal M}_A$ and $d\in C$. Taking $N=A$ and $n=1$,
one obtains
$$1\ot d=\sum_l a_l\theta(c_{l}\ot d_{(1)})\ot d_{(2)}.$$
Applying $\varepsilon_C$ to the second factor, one finds \eqref{2.5.1a}.\\
Conversely, suppose that $\theta\in V_1$
and $z\in W_1$ satisfy \eqref{2.5.1a} and \eqref{2.5.1b}. \eqref{2.5.1b}
implies \eqref{2.5.2}, and \eqref{2.5.1a}
implies \eqref{2.5.3}. Let $\nu=\alpha^{-1}(\theta)$, $\zeta=\beta^{-1}(z)$.
Then (\ref{eq:1a.1.1}-\ref{eq:1a.1.2}) hold, and $(F,G)$ is a Frobenius pair.
\end{proof}

In \cite{CaenepeelMZ97b} it is shown that if
$(H,A,C)$ is a Doi-Hopf structure, $A$ is faithfully flat as a $k$-module,
and $C$ is projective as a $k$-module, then $C$ is finitely generated.
The next proposition shows that, in fact, one does not need the assumption that
$C$ is projective.

\begin{proposition}\prlabel{2.6}
Let $(A,C,\psi)$ be a right-right entwining structure.
If $(F,G)$ is a Frobenius pair,
then $A\ot C$ is finitely generated and projective as a left $A$-module.
\end{proposition}

\begin{proof}
Let $\theta$ and $z=\sum_l a_l\ot c_l$ be as in \thref{2.5}. Then
for all $d\in C$,
\begin{eqnarray*}
1\ot d&=& \psi(d\ot 1)\\
&=& \psi(d_{(1)}\ot \varepsilon(d_{(2)})1)\\
{\rm \eqref{2.5.1b}}~~~~&=&
\sum_l \psi\Bigl(d_{(1)}\ot a_{l\psi}\theta(d_{(2)}^{\psi}\ot c_l)\Bigr)\\
{\rm \eqref{1.1.1}}~~~~&=&
\sum_l a_{l\psi\Psi}\theta(d_{(2)}^{\psi}\ot c_l)_{\Psi'}\ot
d_{(1)}^{\Psi\Psi'}\\
{\rm \eqref{1.1.3}}~~~~&=&
\sum_l a_{l\psi}\theta((d^{\psi})_{(2)}\ot c_l)_{\Psi'}\ot
(d^{\psi})_{(1)}^{\Psi'}\\
{\rm \eqref{2.1.4}}~~~~&=&
\sum_l a_{l\psi}\theta(d^{\psi}\ot c_{l(1)})\ot c_{l(2)}.
\end{eqnarray*}
Write
$c_{l(1)}\ot c_{l(2)}=\sum_{j=1}^{m_l} c_{lj}\ot c'_{lj}$ and
for all $l,j$ consider the map
$$\sigma_{lj}:\ A\ot C\to A,~~~\sigma_{lj}(a\ot d)=
aa_{l\psi}\theta(d^{\psi}\ot c_{lj}).$$
Then for all $a\in A$ and $d\in C$,
$$a\ot d=\sum_{l,j} \sigma_{lj}(a\ot d)(1\ot c'_{lj}),$$
so $\{\sigma_{lj},1\ot c'_{lj}~|~l=1,\cdots,n,~j=1,\cdots,m_l\}$ is
a finite dual basis for $A\ot C$ as a left $A$-module.
\end{proof}

In some situations, one can conclude that $C$ is finitely generated and
projective as a $k$-module.

\begin{corollary}\colabel{2.7}
Let $(A,C,\psi)$ be a right-right entwining structure, and assume that
$(F,G)$ is a Frobenius pair.\\
1) If $A$ is faithfully flat as a $k$-module, then $C$ is finitely generated
as a $k$-module.\\
2) If $A$ is commutative and faithfully flat as a $k$-module, then
$C$ is finitely generated projective as a $k$-module.\\
3) If $k$ is a field, then $C$ is finite dimensional as
a $k$-vector space.\\
4) If $A=k$, then $C$ is finitely generated projective as a $k$-module.
\end{corollary}

\begin{proof}
1) With notation as in \prref{2.6}, let $M$ be the $k$-module generated
by the $c'_{lj}$. Then for all $d\in C$,
$$1\ot d=\sum_{l,j} \sigma_{lj}(1\ot d)\ot c'_{lj}\in A\ot M.$$
Since $A$ is faithfully flat, it follows that $d\in M$, hence $M=C$ is
finitely generated.\\
2) From descent theory: if a $k$-module becomes finitely generated and
projective after a faithfully flat commutative base extension, then it
is itself finitely generated and projective.\\
3) Follows immediately from 1): since $k$ is a field, $A$ is faithfully
flat as a $k$-module, and $C$ is projective as a $k$-module.\\
4) Follows immediately from 2).
\end{proof}

Now we want to recover \cite[Theorem 2.4]{CaenepeelMZ97b} and
\cite[Proposition 3.5]{Brzezinski2}. Assume that $C$ is finitely
generated and projective as a $k$-module, and let
$\{d_i,d_i^*~|~i=1,\cdots,m\}$ be a finite dual basis for $C$.
Then  $C^*\ot A$ can be made into an object of ${}_A{\cal M}(\psi)_A^C$
as follows: for all $a,b,b'\in A,~c^*\in C^*$,
\begin{eqnarray}
b(c^*\ot a)b'&=& \sum_i\lan c^*,d_i^{\psi}\ran d_i^*\ot
b_{\psi}ab',\eqlabel{2.7.6}\\
\rho^r(c^*\ot a)&=& \sum_i d_i^**c^*\ot a_{\psi}\ot d_i^{\psi}.
\eqlabel{2.7.7}
\end{eqnarray}
This can be checked directly. An
explanation for this at first sight artificial structure is given in
\seref{5}. We now give alternative descriptions for $V$ and $W$.
Recall from \cite{CaenepeelIMZ98} that there are many possibilities to
describe $V$. As we have seen, a natural transformation $\nu\in V$ is
completely determined by $\theta$. Nevertheless,  the
maps $\ol{\nu},~\ul{\nu},~\ol{\lambda}$ or $\ul{\lambda}$ (with notation
as in \prref{2.2}) are possible alternatives. The map $\ol{\lambda}:\ C\ot
A\ot C\to A$
induces $\ol{\phi}:\ A\ot C\to C^*\ot A\cong \Hom(C,A)$. This is the
map we need. At some place it is convenient to use $C^*\ot A$
as the image space, at some other we  prefer $\Hom(C,A)$. Note
that $\ol{\phi}$ is given by
$$\ol{\phi}(a\ot c)(d)=\ol{\lambda}(d\ot a\ot c)=
\ul{\lambda}(a_{\psi}\ot d^{\psi}\ot c)=a_{\psi}\theta(d^{\psi}\ot c),$$
or
\begin{equation}\eqlabel{2.8.1}
\ol{\phi}(a\ot c)=\sum_i d_i^*\ot a_{\psi}\theta(d_i^{\psi}\ot c).
\end{equation}
It turns out that $\ol{\phi}$ is a morphism in ${}_A{\cal M}_A^C$. More
specifically, one has

\begin{proposition}\prlabel{2.8}
Let $(A,C,\psi)$ be a right-right entwining structure. If
$C$ is a finitely generated and projective  $k$-module, then
$$V\cong V_1\cong V_2=\Hom_{AA}^{kC}(A\ot C,C^*\ot A).$$
The isomorphism is
$\alpha_1:\ V_1\to V_2$, with $\alpha_1(\theta)=\ol{\phi}$ given by
\eqref{2.8.1}. The inverse of $\alpha_1$ is
\begin{equation}\eqlabel{2.8.2}
\alpha_1^{-1}(\ol{\phi})(d\ot c)=\ol{\phi}(1\ot c)(d).
\end{equation}
\end{proposition}

\begin{proof}
We first show that $\ol{\phi}\in V_2$. For all $a,b\in A$ and $c\in C$,
we have
\begin{eqnarray*}
b\ol{\phi}(a\ot c)&=&
b\Bigl(\sum_i d_i^*\ot a_{\psi}\theta(d_i^{\psi}\ot c)\Bigr)\\
&=& \sum_{i,j} \lan d_i^*,d_j^{\Psi}\ran d_j^*\ot
b_{\Psi}a_{\psi}\theta(d_i^{\psi}\ot c)\\
&=& \sum_j d_j^*\ot b_{\Psi}a_{\psi}\theta(d_j^{\Psi\psi}\ot c)\\
{\rm \eqref{1.1.1}}~~~~&=& \sum_j d_j^*\ot (ba)_{\psi}\theta(d_j^{\psi}\ot c)\\
&=& \ol{\phi}(ba\ot c)=\ol{\phi}(b(a\ot c)),
\end{eqnarray*}
proving that $\ol{\phi}$ is left $A$-linear. It is also right $A$-linear
because
\begin{eqnarray*}
\ol{\phi}(a\ot c)b&=&\sum_i d_i^*\ot a_{\psi}\theta(d_i^{\psi}\ot c)b\\
{\rm \eqref{2.1.3}}~~~~&=& \sum_i d_i^*\ot a_{\psi}b_{\Psi'\Psi}
\theta(d_i^{\psi\Psi}\ot c^{\Psi'})\\
{\rm \eqref{1.1.1}}~~~~&=&\sum_i d_i^*\ot (ab_{\Psi'})_{\psi}
\theta(d_i^{\psi}\ot c^{\Psi'})\\
&=& \ol{\Phi}(ab_{\Psi'}\ot c^{\Psi'})=\ol{\phi}(a\ot cb).
\end{eqnarray*}
 Notice that the dual basis for $C$ satisfies the
following equality (the proof is left to the reader):
\begin{equation}\eqlabel{2.8.3}
\sum_i\Delta(d_i)\ot d_i^*=\sum_{i,j} d_i\ot d_j\ot d_i^**d_j^*.
\end{equation}
Using this equality one computes
\begin{eqnarray*}
\rho^r(\ol{\phi}(a\ot c))&=&
\rho^r\Bigl(\sum_i d_i^*\ot a_{\psi}\theta(d_i^{\psi}\ot c)\Bigr)\\
{\rm \eqref{2.7.7}}~~~~&=&
\sum_{i,j} d_j^**d_i^*\ot \Bigl(a_{\psi}\theta(d_i^{\psi}\ot c)\Bigr)_{\Psi}
\ot d_j^{\Psi}\\
{\rm \eqref{2.8.3}}~~~~&=&
\sum_i d_i^*\ot \Bigl(a_{\psi}\theta(d_{i(2)}^{\psi}\ot c)\Bigr)_{\Psi}
\ot d_{i(1)}^{\Psi}\\
{\rm \eqref{1.1.1}}~~~~&=&\sum_i d_i^*\ot a_{\psi\Psi}
\theta(d_{i(2)}^{\psi}\ot c)_{\Psi'}\ot d_{i(1)}^{\Psi\Psi'}\\
{\rm \eqref{1.1.3}}~~~~&=&\sum_i d_i^*\ot a_{\psi}
\theta((d_{i}^{\psi})_{(2)}\ot c)_{\Psi'}\ot (d_{i}^{\psi})_{(1)}^{\Psi'}\\
{\rm \eqref{2.1.4}}~~~~&=&\sum_i d_i^*\ot a_{\psi}
\theta(d_i^{\psi}\ot c_{(1)})\ot c_{(2)}\\
&=& \ol{\phi}(a\ot c_{(1)})\ot c_{(2)}.
\end{eqnarray*}
This proves that $\ol{\phi}$ is right $C$-colinear.
Conversely, given $\ol{\phi}\in V_2$, first one needs to show that
$\theta=\alpha_1^{-1}(\ol{\phi})\in V_1$. It is now more convenient
to work with $\Hom(C,A)$ rather than $C^*\ot A$. For $f\in \Hom(C,A)$,
$b,b'\in A$, \eqref{2.7.6} can be rewritten as
\begin{equation}\eqlabel{2.8.4}
(bfb')(c)=b_{\psi}f(c^{\psi})b'.
\end{equation}
Take any $c,d\in C$, $a\in A$ and compute
\begin{eqnarray*}
\theta(c\ot d)a&=&\Bigl(\ol{\phi}(1\ot d)(c)\Bigr)a\\
{\rm \eqref{2.8.4}}~~~~&=&(\ol{\phi}(1\ot d)a)(c)\\
(\ol{\phi}~{\rm is~ right~}A{\hbox{-}}{\rm linear})~~&=&
\ol{\phi}(a_{\psi}\ot d^{\psi})(c)\\
(\ol{\phi}~{\rm is~ left~}A{\hbox{-}}{\rm linear})~~&=&
\Bigl(a_{\psi}\ol{\phi}(1\ot d^{\psi})\Bigr)(c)\\
{\rm \eqref{2.8.4}}~~~~&=&
a_{\psi\Psi}\Bigl(\ol{\phi}(1\ot d^{\psi})(c^{\Psi})\Bigr)\\
&=& a_{\psi\Psi}\ol{\phi}(c^{\Psi}\ot d^{\psi}).
\end{eqnarray*}
This proves that $\theta$ satisfies \eqref{2.1.3}. Before proving
\eqref{2.1.4},
we look at the right $C$-coaction $\rho^r$ on $f=c^*\ot a\in \Hom(C,A)\cong
C^*\ot A$.
Write
$\rho^r(f)=f_{[0]}\ot f_{[1]}\in \Hom(C,A)\ot C.$
Using \eqref{2.7.7}, we find, for all $c\in C$,
\begin{eqnarray*}
f_{[0]}(c)\ot f_{[1]} &=&
\sum_i \lan d_i^**c^*,c\ran a_{\psi}d_i^{\psi}\\
&=& \sum_i  \lan d_i^*,c_{(1)}\ran \lan c^*,c_{(2)}\ran a_{\psi}d_i^{\psi}\\
&=& \lan c^*,c_{(2)}\ran a_{\psi}c_{(1)}^{\psi}\\
&=& \psi(c_{(1)}\ot f(c_{(2)})).
\end{eqnarray*}
This means that for all $f\in \Hom(C,A)$
\begin{equation}\eqlabel{2.8.5}
f_{[0]}(c)\ot f_{[1]}=\psi(c_{(1)}\ot f(c_{(2)})).
\end{equation}
This can be used to show that $\theta$ satisfies \eqref{2.1.4}. Explicitly,
\begin{eqnarray*}
\theta(c_{(2)}\ot d)_{\psi}\ot c_{(1)}^{\psi}&=&
\psi(c_{(1)}\ot \theta(c_{(2)}\ot d))\\
{\rm \eqref{2.8.5}}~~~~&=&
\ol{\phi}(1\ot d)_{[0]}(c)\ot \ol{\phi}(1\ot d)_{[1]}\\
(\ol{\phi}~{\rm is~ right~}C{\hbox{-}}{\rm colinear})~~&=&
\ol{\phi}(1\ot d_{(1)})(c)\ot d_{(2)}\\
&=& \theta(c\ot d_{(1)})\ot d_{(2)}.
\end{eqnarray*}
It remains to be shown that $\alpha_1$ and $\alpha_1^{-1}$
are inverses of each other. First take $\theta\in V_1$. Then for all
$c,d\in C$,
\begin{eqnarray*}
\Bigl((\alpha_1^{-1}\circ \alpha_1)(\theta)\Bigr)(d\ot c)&=&
\alpha_1(\theta)(1\ot c)(d)\\
&=& \sum_i \lan d_i^*,d\ran 1_{\psi}\theta(d_i^{\psi}\ot c)\\
&=& \theta(d\ot c).
\end{eqnarray*}
Finally, for $\ol{\phi}\in V_2$, $a\in A$ and $c,d\in C$:
\begin{eqnarray*}
\Bigl((\alpha_1\circ \alpha_1^{-1})(\ol{\phi})\Bigr)(a\ot c)(d)
&=& \sum_i \lan d_i^*,d\ran a_{\psi}\alpha_1^{-1}(\ol{\phi})(d_i^{\psi}\ot c)\\
&=& \sum_i \lan d_i^*,d\ran a_{\psi}\ol{\phi}(1\ot c)(d_i^{\psi})\\
&=& a_{\psi}\ol{\phi}(1\ot c)(d^{\psi})\\
&=& (a\ol{\phi}(1\ot c))(d)=\ol{\phi}(a\ot c)(d).
\end{eqnarray*}
\end{proof}

Now we give an alternative description for $W_2$.

\begin{proposition}\prlabel{2.9}
Let $C$ be finitely generated and projective as
a $k$-module. Then
$$W\cong W_1\cong W_2=\Hom_{AA}^{kC}(C^*\ot A,A\ot C).$$
The isomorphism $\beta_1:\ W_1\to W_2$ is given by $\beta_1(z)=\phi$
with
$$\phi(c^*\ot a)=\sum_l a_la_{\psi}\ot
\lan c^*,c_{l(2)}\ran c_{l(1)}^{\psi},$$
and the inverse of $\beta_1$ is given by
\begin{equation}\eqlabel{2.9.0}
\beta_1^{-1}(\phi)=\phi(\varepsilon\ot 1).
\end{equation}
\end{proposition}

\begin{proof}
We have to show that $\beta_1(z)=\phi$ is left and right $A$-linear
and right $C$-colinear.
For all $c^*\in C^*$ and $a,b\in A$,
\begin{eqnarray*}
\Phi(c^*\ot ab)&=& \sum_l a_l(ab)_{\psi}\ot \lan c^*,c_{j(2)}\ran
c_{j(1)}^{\psi}\\
{\rm \eqref{1.1.1}}~~~~&=&
\sum_l a_la_{\psi}b_{\Psi}\ot \lan c^*,c_{j(2)}\ran c_{j(1)}^{\psi\Psi}\\
&=& \Bigl(\sum_l a_la_{\psi}\ot \lan c^*,c_{j(2)}\ran c_{j(1)}^{\psi}\Bigr)b\\
&=& \Phi(c^*\ot a)b,
\end{eqnarray*}
proving that $\Phi$ is right $A$-linear. The proof of left $A$-linearity
goes as follows:
\begin{eqnarray*}
\Phi(b(c^*\ot a))&=&
\sum_i \Phi\Bigl(\lan c^*,d_i^{\psi}\ran d_i^*\ot b_{\psi}a\Bigr)\\
&=& \sum_{i,l}\lan c^*,d_i^{\psi}\ran a_l(b_{\psi}a)_{\Psi}\ot
\lan d_i^*,c_{l(2)}\ran c_{l(1)}^{\Psi}\\
{\rm \eqref{1.1.1}}~~~~&=&
\sum_{i,l}\lan c^*,d_i^{\psi}\ran a_lb_{\psi\Psi}a_{\Psi'}\ot
\lan d_i^*,c_{l(2)}\ran c_{l(1)}^{\Psi\Psi'}\\
&=&
\sum_{l}\lan c^*,c_{l(2)}^{\psi}\ran a_lb_{\psi\Psi}a_{\Psi'}\ot
c_{l(1)}^{\Psi\Psi'}\\
{\rm \eqref{1.1.3}}~~~~&=&
\sum_{l}\lan c^*,(c_{l}^{\psi})_{(2)}\ran a_lb_{\psi}a_{\Psi'}\ot
(c_{l}^{\psi})_{(1)}^{\Psi'}\\
{\rm \eqref{2.3.1}}~~~~&=&
\sum_{l}\lan c^*,(c_{l(2)}\ran ba_la_{\Psi'}\ot
(c_{l(1)}^{\Psi'}\\
&=& b\Phi(c^*\ot a).
\end{eqnarray*}
Next one needs to show that $\Phi$ is right $C$-colinear. Using \eqref{2.8.5},
one finds
\begin{eqnarray*}
\ol{\phi}((c^*\ot a)_{[0]})\ot (c^*\ot a)_{[1]}
&=& \sum_i \ol{\phi}(d_i^**c^*\ot a_{\psi})\ot d_i^{\psi}\\
&=& \sum_{i,l} a_la_{\psi\Psi}\ot\lan d_i^**c^*,c_{l(2)}\ran
c_{l(1)}^{\Psi}\ot d_i^{\psi}\\
&=& \sum_{i,l} a_la_{\psi\Psi}\ot\lan d_i^*,c_{l(2)}\ran\lan c^*,c_{l(3)}\ran
c_{l(1)}^{\Psi}\ot d_i^{\psi}\\
&=& \sum_{l} a_la_{\psi\Psi}\ot\lan c^*,c_{l(3)}\ran
c_{l(1)}^{\Psi}\ot c_{l(2)}^{\psi}\\
{\rm \eqref{1.1.3}}~~~~&=&
\sum_{l} a_la_{\psi}\ot\lan c^*,c_{l(2)}\ran (c_l^{\psi})_{(1)}
\ot (c_l^{\psi})_{(2)}\\
&=& \rho^r(\ol{\phi}(c^*\ot a)).
\end{eqnarray*}
Conversely, let $\phi\in W_2$ and put $z=\phi(\varepsilon\ot 1)=
\sum_l a_l\ot c_l$. Using \eqref{2.7.6}, we see that
$a(\varepsilon\ot 1)=(\varepsilon\ot 1)a$, for all $a\in A$, hence
$az=a\phi(\varepsilon\ot 1)=\phi(a(\varepsilon\ot 1))
=\phi((\varepsilon\ot 1)a)=\phi(\varepsilon\ot 1)a=za$, and $z\in W_1$.\\
Take $z=\sum_l a_l\ot c_l\in W_1$. Then
$$\beta_1^{-1}(\beta_1(z))=\sum_l a_l1_{\psi}\ot \lan\varepsilon,c_{l(2)}\ran
c_{l(1)}^{\psi}=z.$$
Finally, take $\phi\in W_2$, and write
$\beta_1^{-1}(\phi)=\phi(\varepsilon\ot 1)=\sum_l a_l\ot c_l$.
$C^*\ot A$ and $A\ot C$ are right $C$-comodules and left $C^*$-modules.
Since $\phi$ is right $A$-linear, right $C$-colinear and left
$C^*$-linear,
\begin{eqnarray*}
&&\hspace*{-12mm}
\phi(c^*\ot a)=\phi(c^*\ot 1)a=(c^*\cdot\phi(\varepsilon\ot 1))a=
(c^*\cdot (\sum_l a_l\ot c_l))a\\
&=& \lan c^*,c_{l(2)}\ran \sum a_la_{\psi}\ot c_{l(1)}^{\psi}=
\beta_1(z)(c^*\ot a)
\end{eqnarray*}
and it follows that $\phi=\beta_1(z)=\beta_1(\beta_1^{-1}(\phi))$,
as required.
\end{proof}

Suppose that $C$ is finitely generated and projective as a $k$-module.
>From \prref{2.8}, it follows that $F$ is separable if and only if there
exists a map $\ol{\phi}\in V_2=\Hom_{A,A}^C(A\ot C,C^*\ot A)$ such that
$\ol{\phi}(1\ot c_{(2)})(c_{(1)})=\varepsilon(c)1$, for all $c\in C$.
In the Doi-Hopf case, this implies the Maschke Theorem
in \cite{CaenepeelMZ97a}. Now we apply the same procedure to determine
when $(F,G)$ is  a Frobenius pair.

\begin{theorem}\thlabel{2.10}
Consider an entwining structure $(A,C,\psi)$, and assume that $C$ is
finitely generated projective as a $k$-module.
Let $F:\ {\cal M}(\psi)_A^C\to {\cal M}_A$ be the functor forgetting
the $C$-coaction, and $G=\bullet\ot C$ be its right adjoint. Then the
following statements are equivalent:\\
1) $(F,G)$ is a Frobenius pair.\\
2) There exist $z=\sum a_l\ot c_l\in W_1$ and $\theta\in V_1$ such that the
maps
$$\phi:\ C^*\ot A\to A\ot C~~{\rm and}~~\ol{\phi}:\
A\ot C\to C^*\ot A,$$
given by
\begin{eqnarray}
\phi(c^*\ot a)&=& \sum_l a_la_{\psi}\ot \lan
c^*,c_{l(2)}\ran c_{l(1)}^{\psi},
\eqlabel{2.10.1}\\
\ol{\phi}(a\ot c)&=& \sum_i d_i^*\ot a_{\psi}\theta(d_i^{\psi}\ot c),
\eqlabel{2.10.2}
\end{eqnarray}
are inverses of each other.\\
3) $C^*\ot A$ and $A\ot C$ are isomorphic as objects in
${}_A{\cal M}(\psi)_A^C$.
\end{theorem}

\begin{proof}
$\ul{1)\Rightarrow 2)}$. Let $z\in W_1$ and
$\theta\in V_1$ be as in \thref{2.5}. Then $\phi=\beta_1(z)$ and
$\ol{\phi}=\alpha_1(\theta)$ are morphisms in ${}_A{\cal M}(\psi)_A^C$,
and
\begin{eqnarray*}
\ol{\phi}(\phi(\varepsilon\ot 1))=\ol{\phi}(z)&=&
\sum_{i,l} d_i^*\ot a_{l\psi}\theta(d_i^{\psi}\ot c_l)\\
{\rm \eqref{2.5.1b}}~~~~&=&
\sum_{i,l} d_i^*\ot \varepsilon(d_i)1=\varepsilon\ot 1.
\end{eqnarray*}
The fact that $\phi$ and $\ol{\phi}$ are right $A$-linear and
left $C^*$-linear implies that $\ol{\phi}\circ\phi=I_{C^*\ot A}$.
Similarly,  for all $c\in C$,
\begin{eqnarray*}
\phi(\ol{\phi}(1\ot c))&=&
\phi(\sum_i d_i^*\ot \theta(d_i\ot c))\\
&=& \sum_{i,l}a_l\theta(d_i\ot c)_{\psi}\ot \lan d_i^*,c_{l(2)}\ran
c_{l(1)}^{\psi}\\
&=& \sum_{l}a_l\theta(c_{l(2)}\ot c)_{\psi}\ot
c_{l(1)}^{\psi}\\
{\rm \eqref{2.1.4}}~~~~&=& \sum_{l}a_l\theta(c\ot c_{(1)})\ot
c_{(2)}\\
{\rm \eqref{2.5.1a}}~~~~&=& \varepsilon(c_{(1)})1\ot c_{(2)})\\
&=& 1\ot c.
\end{eqnarray*}
Since $\phi$ and $\ol{\phi}$ are left $A$-linear,
$\phi\circ\ol{\phi}=I_{A\ot C}$.\\
$\ul{2)\Rightarrow 3)}$. Obvious, since $\phi$ and $\ol{\phi}$ are
in ${}_A{\cal M}(\psi)_A^C$.\\
$\ul{3)\Rightarrow 1)}$. Let $\phi:\ C^*\ot A\to A\ot C$ be the connecting
isomorphism, and put $z=\phi(\varepsilon\ot 1)=\sum_l a_l\ot c_l\in W_1$,
$\theta=\alpha_1^{-1}(\phi^{-1})\in V_1$. Applying \eqref{2.8.1}
and \eqref{2.9.0}, one finds
$$\varepsilon\ot 1=\phi^{-1}(\phi(\varepsilon\ot 1))=
\sum_i d_i^*\ot a_{l\psi}\theta(d_i^{\psi}\ot c_l).$$
Evaluating this equality at $c\in C$, one obtains \eqref{2.5.1b}.
For all $c\in C$,
$$1\ot c=\phi(\phi^{-1}(1\ot c))=\sum_l
a_l\theta(c_l\ot c_{(1)})\ot c_{(2)}.$$
Applying $\varepsilon$ to the second factor, one finds \eqref{2.5.1a}.
\thref{2.5} implies that $(F,G)$ is a Frobenius pair.
\end{proof}

\begin{remark}\rm
Recently M.\ Takeuchi observed that entwined modules can be viewed as
comodules over certain corings. This observation has been exploited in
 \cite{Brzezinski00} to derive some properties of coring counterparts
of functors $F$ and $G$.
It is quite clear  that the procedure applied in \seref{1a} to extension
and restriction of scalars  can be adapted to functors associated to
corings, leading to a generalization of the results in this Section.
This will be the subject of a future publication.
\end{remark}

\section{The functor forgetting the $A$-action}\selabel{4}
Again, let $(A,C,\psi)$ be a right-right entwining structure. The functor
$G':\ {\cal M}(\psi)_A^C\to {\cal M}^C$ forgetting the $A$-action has
a left adjoint $F'$. The unit $\mu$ and the counit $\eta$ of the adjunction
are given at the end of \seref{1}.

\begin{lemma}\lelabel{4.1}
Let $M\in {}_A{\cal M}(\psi)_A^C$, $N\in {}^C{\cal M}(\psi)_A^C$.
Then $F'G'(M)\in {}_A{\cal M}(\psi)_A^C$ and
$G'F'\in {}^C{\cal M}(\psi)_A^C$. The left structures are given by
$$a(m\ot b)=am\ot b~~{\rm and}~~\rho^l(n\ot b)=
\sum n_{[-1]}\ot n_{[0]}\ot b,$$
for all $a,b\in A$, $m\in M$, $n\in N$. Furthermore $\mu_M$ is left
$A$-linear, and $\nu_N$ is left $C$-colinear.
\end{lemma}

Now write $V'=\dul{\rm Nat}(G'F',1_{{\cal M}^C})$,
$W'=\dul{\rm Nat}(1_{{\cal C}}, F'G')$. Following the philosophy of the
previous Sections, we give more explicit descriptions of $V'$ and $W'$.
We do not give detailed  proofs, however, since the arguments are dual to the
ones in the previous Section. Let
\begin{equation}\eqlabel{4.2.0}
V'_1=\{\vartheta\in (C\ot A)^*~|~
\vartheta(c_{(1)}\ot a_{\psi})c_{(2)}^{\psi}=
\vartheta(c_{(2)}\ot a)c_{(1)},~{\rm for~all~}c\in C,a\in A\}.
\end{equation}

\begin{proposition}\prlabel{4.2}
The map $\alpha:\ V'\to V'_1$, $\alpha(\nu')=\varepsilon\circ \nu_C$
is an isomorphism.
\end{proposition}

\begin{proof}
Details are left to the reader. Given $\vartheta\in
V'$,  for $N\in {\cal M}^C$, the natural map
$\nu'_N:\ N\ot A\to N$ is
$$\nu'_N(n\ot a)=\sum \vartheta(n_{[1]}\ot a)n_{[0]}$$
\end{proof}

For any $k$-linear map $e:\ C\to A\ot A$ and $c\in C$, we use the notation
$e(c)=e^1(c)\ot e^2(c)$ (summation understood).
Let $W'_1$ be the $k$-submodule of $\Hom(C,A\ot A)$ consisting of
maps $e$ satisfying
\begin{eqnarray}
e^1(c_{(1)})\ot e^2(c_{(1)})\ot c_{(2)}&=&
e^1(c_{(2)})_{\psi}\ot e^2(c_{(2)})_{\Psi}\ot c_{(1)}^{\psi\Psi},
\eqlabel{4.3.1}\\
e^1(c)\ot e^2(c)a&=& a_{\psi}e^1(c^{\psi})\ot e^2(c^{\psi}).\eqlabel{4.3.2}
\end{eqnarray}

\begin{proposition}\prlabel{4.3}
The map $\beta:\ W'\to W'_1$ given by
\begin{eqnarray*}
\beta(\zeta')&=& (\varepsilon\ot I_A\ot I_A)\circ \zeta'_{A\ot C}
\circ (\eta_A\ot I_C)\\
&=& (I_A\ot \varepsilon\ot I_A)\circ \zeta'_{C\ot A}
\circ (I_C\ot \eta_A)
\end{eqnarray*}
is an isomorphism. Given $e\in W'_1$, $\zeta'=\beta^{-1}(e)$ is
recovered from $e$ as follows: for $M\in {\cal M}(\psi)_A^C$,
$$\zeta'_M(m)= \sum m_{[0]}e^1(m_{[1]})\ot e^2(m_{[1]}).$$
\end{proposition}

\begin{proof}
We show that $\beta$ is well-defined, leaving other details to
the reader. Consider a commutative diagram
$$\begin{diagram}
C&\rTo^{I_C\ot \eta_A}&C\ot A&\rTo^{\zeta'_{C\ot A}}&C\ot A\ot A\\
\dTo^{I_C}&&\dTo^{\psi}&&\dTo^{\psi\ot I_A}\\
C&\rTo^{\eta_A\ot I_C}&A\ot C&\rTo^{\zeta'_{A\ot C}}&A\ot C\ot A
\end{diagram}$$
The map $\ol{\lambda}=\zeta'_{C\ot A}\circ (I_C\ot \eta_A)$ is left and right
$C$-colinear. Write
$\ol{\lambda}(c)=\sum_i c_i\ot a_i\ot a'_i$.
Then
$$c_{(1)}\ot \ol{\lambda}(c_{(2)})=
\sum_i c_{i(1)}\ot c_{i(2)}\ot a_i\ot a'_i.$$
Applying $\varepsilon$ to the second factor, one finds
$$c_{(1)}\ot e(c_{(2)})= \ol{\lambda}(c).$$
The right $C$-colinearity of $\ol{\lambda}$ implies that
\begin{eqnarray*}
\ol{\lambda}(c_{(1)})\ot c_{(2)}&=&
\sum_i a_{i\psi}\ot a'_{i\Psi}\ot c_i^{\psi\Psi}\\
&=& e^1(c_{(2)})_{\psi}\ot e^2(c_{(2)})_{\Psi}\ot c_{(1)}^{\psi\Psi},
\end{eqnarray*}
and hence proves \eqref{4.3.1}. To prove \eqref{4.3.2} note that
$\ul{\lambda}=\zeta'_{A\ot C}\circ (\eta_A\ot I_C)$ is left and right
$A$-linear, hence
\begin{eqnarray*}
e^1(c)\ot e^2(c)a&=&
(I_A\ot\varepsilon\ot I_A)(\zeta'_{A\ot C}((1\ot c)a))\\
&=&(I_A\ot\varepsilon\ot I_A)(\zeta'_{A\ot C}(a_{\psi}\ot c^{\psi})\\
&=&a_{\psi}(I_A\ot\varepsilon\ot I_A)(\zeta'_{A\ot C}(1\ot c^{\psi})\\
&=&a_{\psi}e^1(c^{\psi})\ot e^2(c^{\psi}).
\end{eqnarray*}
\end{proof}

\begin{proposition}\prlabel{4.4}
Let $(A,C,\psi)$ be a right-right entwining structure.\\
1) $F'=\bullet\ot A:\ {\cal M}^C\to {\cal M}(\psi)_A^C$ is separable if and
only if
there exists $\vartheta\in V'_1$ such that for all $c\in C$,
\begin{equation}\eqlabel{4.4.1}
\vartheta(c\ot 1)=\varepsilon(c).
\end{equation}
2) $G':\ {\cal M}(\psi)_A^C\to {\cal M}^C$ is separable if and only if
there exists $e\in W'_1$ such that for all $c\in C$,
\begin{equation}\eqlabel{4.4.2}
e^1(c)e^2(c)=\varepsilon(c)1.
\end{equation}
3) $(F',G')$ is a Frobenius pair if and only if there exist
$\vartheta\in V'_1$ and $e\in W'_1$ such that
\begin{eqnarray}
\varepsilon(c) 1&=& \vartheta(c_{(1)}\ot e^1(c_{(2)})) e^2(c_{(2)})
\eqlabel{4.4.3}\\
&=& \vartheta(c_{(1)}^{\psi}\ot e^2(c_{(2)})) e^1(c_{(2)})_{\psi}.
\eqlabel{4.4.4}
\end{eqnarray}
\end{proposition}

\begin{proof}
We only prove 3). If $(F',G')$ is a Frobenius pair, then there exist
$\nu'\in V'$ and $\zeta'\in W'$ such that \eqref{1a.1.1} and \eqref{1a.1.2}
hold. Take $\vartheta\in V'_1$ and $e\in W'_1$ corresponding to
$\nu'$ and $\zeta'$ and write down \eqref{1a.1.1} applied to $n\ot 1$
with $n\in N\in {\cal M}^C$,
\begin{equation}\eqlabel{4.4.5}
n\ot 1=\bigl((\nu'_N\ot I_A)\circ \zeta'_{N\ot A})\bigr)(n\ot 1)=
\vartheta(n_{[1]}\ot e^1(n_{[2]}))n_{[0]}\ot e^1(n_{[2]})).
\end{equation}
Taking $N=C$, $n=c$, and applying $\varepsilon_C$ to the first factor,
one obtains \eqref{4.4.3}. Conversely, if
$\vartheta\in V'_1$ and $e\in W'_1$ satisfy \eqref{4.4.3}, then
\eqref{4.4.5} is satisfied for all $N\in {\cal M}^C$, and \eqref{1a.1.1}
follows since
$\nu'_N\ot I_A$ and $\zeta'_{N\ot A}$ are right $A$-linear.\\
Now write down \eqref{1a.1.2} applied to $m\in M\in {\cal M}(\psi)_A^C$,
\begin{equation}\eqlabel{4.4.6}
m=\bigl(\nu'_{G'(M\ot A)}\circ G'(\zeta'_M)\bigr)(m)=
\theta(m_{[1]}^{\psi}\ot e^2(m_{[2]}))m_{[0]}e^1(m_{[2]})_{\psi}.
\end{equation}
Take $M=C\ot A$, $m=c\ot 1$, and apply $\varepsilon_C$ to the first factor.
This gives \eqref{4.4.4}. Conversely,
if  $\vartheta\in V'_1$ and $e\in W'_1$ satisfy \eqref{4.4.4}, then
application of \eqref{4.4.4} to the second and third factors in
$m_{[0]}\ot m_{[1]}\ot 1$, and then $\varepsilon_C$ to the second
factor shows that \eqref{4.4.6} holds for all $M\in {\cal M}(\psi)_A^C$.
Finally note that \eqref{4.4.6} is equivalent to \eqref{1a.1.2}.
\end{proof}

Inspired by the results in the previous Section, we ask the following
question: assuming $(F',G')$ is a Frobenius pair, when is $A$ finitely
generated projective as a $k$-module. We give a partial answer in the
next Proposition. We assume that $\psi$ is bijective (cf.\
\cite[Section~6]{Brzezinski1}). In the Doi-Hopf
case, this is true if the underlying Hopf algebra $H$ has a twisted
antipode. The inverse of $\psi$ is then given by the formula
$$\psi^{-1}(a\ot c)= c\ol{S}(a_{(1)})\ot a_{(0)}.$$

\begin{proposition}\prlabel{4.5}
Let $(A,C,\psi)$ be a right-right entwining structure. With notation
as above, assume that $(F',G')$ is a Frobenius pair. If there exists
$c\in C$ such that $\varepsilon(c)=1$, and if $\psi$ is invertible,
with inverse $\varphi=\psi^{-1}:\ A\ot C\to C\ot A$, then $A$ is
finitely generated and projective as a $k$-module.
\end{proposition}

\begin{proof}
Observe first that $(A,C,\varphi)$ is a left-left entwining structure.
This means that (\ref{eq:1.1.1}-\ref{eq:1.1.4}) hold, but with $A$ and
$C$ replaced by $A^{\rm op}$ and $C^{\rm cop}$. In particular,
\begin{eqnarray}
\varepsilon(c^{\varphi})a_{\varphi}&=&\varepsilon(c)a,\eqlabel{4.5.1}\\
a_{\varphi}\ot \Delta(c^{\varphi})&=&
a_{\varphi\phi}\ot c_{(1)}^{\varphi}\ot c_{(2)}^{\phi}.\eqlabel{4.5.2}
\end{eqnarray}
This can be seen as follows: rewrite \eqref{1.1.2} and \eqref{1.1.3}
as commutative diagrams, reverse the arrows, and replace $\psi$ by
$\varphi$. Then we have \eqref{4.5.1} and \eqref{4.5.2} in diagram form.
Now fix $c\in C$ such that $\varepsilon(c)=1$. Then for all $a\in A$,
\begin{eqnarray*}
a=\varepsilon(c)a&=&\varepsilon(c^{\varphi})a_{\varphi}\\
{\rm \eqref{4.4.3}}~~~~&=&
\vartheta\bigl((c^{\varphi})_{(1)}\ot e^1((c^{\varphi})_{(2)})\bigr)
e^2((c^{\varphi})_{(2)})a_{\varphi}\\
{\rm \eqref{4.5.2}}~~~~&=&
\vartheta\bigl(c^{\varphi}_{(1)}\ot e^1(c^{\phi}_{(2)})\bigr)
e^2(c^{\phi}_{(2)})a_{\varphi\phi}\\
{\rm \eqref{4.3.2}}~~~~&=&
\vartheta\bigl(c^{\varphi}_{(1)}\ot a_{\varphi\phi\psi}
e^1(c^{\phi\psi}_{(2)})\bigr)e^2(c^{\phi\psi}_{(2)})\\
(\varphi=\psi^{-1})~~~~&=&
\vartheta\bigl(c^{\varphi}_{(1)}\ot a_{\varphi}
e^1(c_{(2)})\bigr)e^2(c_{(2)}).\\
\end{eqnarray*}
Write
$(I\ot e)\Delta(c)=\sum_{i=1}^m c_i\ot b_i\ot a_i\in C\ot A\ot A$.
For $i=1,\cdots,m$, define $a_i^*\in A^*$ by
$$\lan a_i^*,a\ran = \vartheta\bigl(c^{\varphi}_i\ot
a_{\varphi}b_i\bigr).$$
Then $\{a_i,a_i^*~|~i=1,\cdots, m\}$ is a finite dual basis of $A$ as
a $k$-module.
\end{proof}

{}From now on we assume that $A$ is finitely generated and projective
with finite dual basis $\{a_i,a_i^*~|~i=1,\cdots, m\}$. The proof of
the next Lemma is straightforward, and therefore left to the reader.

\begin{lemma}\lelabel{4.6}
Let $(A,C,\psi)$ be a right-right entwining structure, and assume that
$A$ is finitely generated and projective as a $k$-module. Then
$A^*\ot C\in {}^C{\cal M}(\psi)_A^C$. The structure is given by the
formulae
\begin{eqnarray}
(a^*\ot c)b&=& \lan a^*, b_{\psi}a_i\ran a_i^*\ot c^{\psi},\eqlabel{4.6.1}\\
\rho^r(a^*\ot c)&=& a^*\ot c_{(1)}\ot c_{(2)},\eqlabel{4.6.2}\\
\rho^l(a^*\ot c)&=& \lan a^*,a_{i\psi}\ran c_{(1)}^{\psi}\ot a_i^*\ot
c_{(2)}^{\psi}.\eqlabel{4.6.3}
\end{eqnarray}
\end{lemma}

We now give alternative descriptions of $V'$ and $W'$.

\begin{proposition}\prlabel{4.7}
Let $(A,C,\psi)$ be a right-right entwining structure, and assume that
$A$ is finitely generated and projective as a $k$-module. Then there is
an isomorphism
$$\beta_1:\ W'_1\to W'_2=\Hom^{CC}_{kA}(A^*\ot C, C\ot A),$$
$\beta_1(e)=\Omega$, with
$$\Omega(a^*\ot c)=\lan a^*,e^1(c_{(2)})_{\psi}\ran c_{(1)}^{\psi}\ot
e^2(c_{(2)}).$$
The inverse of $\beta_1$ is given by $\beta_1^{-1}(\Omega)=e$ with
$$e(c)=\sum_i a_i \ot (\varepsilon_C\ot I_A)\Omega(a^*\ot c).$$
\end{proposition}

\begin{proof}
We first prove that $\beta_1$ is well-defined.\\
a) $\beta_1(e)=\Omega$ is right $A$-linear: for all $a^*\in A^*$,
$c\in C$ and $b\in A$, we have
\begin{eqnarray*}
\Omega\bigl((a^*\ot c)b\bigr)&=&
\sum_i \lan a^*,b_{\psi}a_i\ran \Omega(a_i^*\ot c^{\psi})\\
&=& \sum_i \lan a^*,b_{\psi}a_i\ran \lan a_i^*,
e^1((c^{\psi})_{(2)})_{\Psi}\ran
(c^{\psi})_{(1)}^{\Psi}\ot e^2((c^{\psi})_{(2)})\\
&=& \lan a^*,b_{\psi}e^1((c^{\psi})_{(2)})_{\Psi}\ran
(c^{\psi})_{(1)}^{\Psi}\ot e^2((c^{\psi})_{(2)})\\
{\rm \eqref{1.1.3}}~~~~&=&
\lan a^*,b_{\psi\psi'}e^1(c^{\psi}_{(2)})_{\Psi}\ran
c^{\psi'\Psi}_{(1)}\ot e^2(c^{\psi}_{(2)})\\
{\rm \eqref{1.1.1}}~~~~&=&
\lan a^*,\bigl(b_{\psi}e^1(c^{\psi}_{(2)})\bigr)_{\Psi}\ran
c^{\Psi}_{(1)}\ot e^2(c^{\psi}_{(2)})\\
{\rm \eqref{4.3.2}}~~~~&=&
\lan a^*,e^1(c_{(2)})_{\Psi}\ran
c^{\Psi}_{(1)}\ot e^2(c_{(2)})b\\
&=& \Omega(a^*\ot c)b.
\end{eqnarray*}
b) $\beta_1(e)=\Omega$ is right $C$-colinear: for all $a^*\in A^*$ and
$c\in C$, we have
\begin{eqnarray*}
\rho^r(\Omega(a^*\ot c))&=&
\lan a^*,e^1(c_{(2)})_{\psi}\ran \rho^r(c_{(1)}^{\psi}\ot e^2(c_{(2)})\\
&=&
\lan a^*,e^1(c_{(2)})_{\psi}\ran \bigl(c_{(1)}^{\psi}\bigr)_{(1)}\ot
e^2(c_{(2)})_{\Psi}\ot \bigl(c_{(1)}^{\psi}\bigr)_{(2)}^{\Psi}\\
{\rm \eqref{1.1.3}}~~~~&=&
\lan a^*,e^1(c_{(3)})_{\psi\psi'}\ran c_{(1)}^{\psi'}\ot
e^2(c_{(3)})_{\Psi}\ot c_{(2)}^{\psi\Psi}\\
{\rm \eqref{4.3.1}}~~~~&=&
\lan a^*,e^1(c_{(2)})_{\psi'}\ran c_{(1)}^{\psi'}\ot
e^2(c_{(2)})\ot c_{(3)}\\
&=& \Omega(a^*\ot c_{(1)})\ot c_{(2)}.
\end{eqnarray*}
c) $\beta_1(e)=\Omega$ is left $C$-colinear: for all $a^*\in A^*$ and
$c\in C$, we have
\begin{eqnarray*}
\rho^l(\Omega(a^*\ot c))&=&
\lan a^*, e^1(c_{(2)})_{\psi}\ran (c_{(1)}^{\psi})_{(1)}\ot
(c_{(1)}^{\psi})_{(2)}\ot e^2(c_{(2)})\\
{\rm \eqref{1.1.3}}~~~~&=&
\lan a^*, e^1(c_{(3)})_{\psi\psi'}\ran c_{(1)}^{\psi'}\ot
c_{(2)}^{\psi}\ot e^2(c_{(3)})\\
&=&
\sum_i \lan a^*, a_{i\psi}\ran \lan a_i^*,
e^1(c_{(3)})_{\psi'}\ran c_{(1)}^{\psi'}\ot
c_{(2)}^{\psi}\ot e^2(c_{(3)})\\
&=& \sum_i \lan a^*, a_{i\psi}\ran c_{(1)}^{\psi}\ot
\Omega(a_i^*\ot c_{(2)}).
\end{eqnarray*}
The proof that $\beta_1^{-1}(\Omega)=e$
satisfies \eqref{4.3.1} and \eqref{4.3.2} is left to the reader.
The maps $\beta_1$ and $\beta_1^{-1}$ are inverses of each other since
\begin{eqnarray*}
\beta^{-1}(\beta(e))(c)&=&
\sum_i a_i\ot (\varepsilon_C\ot I_A)
\lan a_i^*,e^1(c_{(2)})_{\psi}\ran c_{(1)}^{\psi}\ot e^2(c_{(2)})\\
&=&\sum_i a_i\ot \lan a_i^*,e^1(c_{(2)})_{\psi}\ran
\varepsilon_C(c_{(1)}^{\psi})\ot e^2(c_{(2)})\\
&=& \lan a_i^*,e^1(c)\ran a_i\ot e^2(c)=e^1(c)\ot e^2(c),
\end{eqnarray*}
\begin{eqnarray*}
\beta(\beta^{-1}(\omega))(a^*\ot c)&=&
\lan a^*,(a_i)_{\psi}\ran c_{(1)}^{\psi}\ot
(\varepsilon_C\ot I_A)\Omega(a_i^*\ot c_{(2)})\\
&=& (I_C\ot\varepsilon_C\ot I_A)\bigl(
\lan a^*,(a_i)_{\psi}\ran c_{(1)}^{\psi}\ot\Omega(a_i^*\ot c_{(2)})\bigr)\\
&=& (I_C\ot\varepsilon_C\ot I_A)\rho^l(\Omega(a^*\ot c))\\
&=& \Omega(a^*\ot c).
\end{eqnarray*}
At the last step, we used that for all $c\in C$ and $a\in A$,
$$(I_C\ot\varepsilon_C\ot I_A)\rho^l(c\ot a)=c\ot a.$$
\end{proof}

\begin{proposition}\prlabel{4.8}
Let $(A,C,\psi)$ be a right-right entwining structure If
$A$ is a finitely generated projective $k$-module, then the map
$$\alpha_1:\ V'_1\to V'_2=\Hom^{CC}_{kA}(C\ot A,A^*\ot C),$$
defined by $\alpha_1(\vartheta)=\ol{\Omega}$, with
$$\ol{\Omega}(c\ot a)=\lan\vartheta,c_{(1)}\ot a_{\psi}a_i\ran a_i^*\ot
c_{(2)}^{\psi}$$
is an isomorphism. The inverse of $\alpha_1$ is given by
$\alpha_1^{-1}(\ol{\Omega})=\vartheta$ with
$$\vartheta(c\ot a)=\lan\ol{\Omega}(c\ot a),1_A\ot \varepsilon_C\ran.$$
\end{proposition}

\begin{proof}
We first show that $\alpha_1$ is well-defined. Take $\vartheta\in V'_1$,
and let $\alpha_1(\vartheta)=\ol{\Omega}$.\\
a) $\ol{\Omega}$ is right $A$-linear since for all $a,b\in A$ and $c\in C$,
\begin{eqnarray*}
\ol{\Omega}(c\ot a)b&=&
\sum_{i,j} \lan\vartheta,c_{(1)}\ot a_{\psi}a_i\ran \lan a_i^*,
b_{\Psi}a_j\ran a_j^*\ot c_{(2)}^{\psi\Psi}\\
&=&\sum_{j} \lan\vartheta,c_{(1)}\ot a_{\psi}
b_{\Psi}a_j\ran a_j^*\ot c_{(2)}^{\psi\Psi}\\
{\rm \eqref{1.1.1}}~~~~&=&
\sum_{j} \lan\vartheta,c_{(1)}\ot
(ab)_{\psi}a_j\ran a_j^*\ot c_{(2)}^{\psi}\\
&=& \ol{\Omega}(c\ot ab).
\end{eqnarray*}
b) $\ol{\Omega}$ is right $C$-colinear since for all $a\in A$ and $c\in C$,
\begin{eqnarray*}
\rho^r(\ol{\Omega}(c\ot a))&=&
\vartheta(c_{(1)}\ot a_{\psi}a_i)a_i^*\ot (c_{(2)}^{\psi})_{(1)}
\ot (c_{(2)}^{\psi})_{(2)}\\
{\rm \eqref{1.1.3}}~~~~&=&
\vartheta(c_{(1)}\ot a_{\psi\Psi}a_i)a_i^*\ot c_{(2)}^{\Psi}
\ot c_{(3)}^{\psi}\\
&=& \ol{\Omega}(c_{(1)}\ot a_{\psi})\ot c_{(2)}^{\psi}.
\end{eqnarray*}
c) $\ol{\Omega}$ is left $C$-colinear since for all $a\in A$ and $c\in C$,
\begin{eqnarray*}
\rho^l(\ol{\Omega}(c\ot a))&=&
\sum_{i,j} \vartheta(c_{(1)}\ot a_{\psi}a_i) \lan a_i^*,a_{j\Psi}\ran
(c_{(2)}^{\psi})_{(1)}^{\Psi}\ot a_j^*\ot (c_{(2)}^{\psi})_{(2)}\\
{\rm \eqref{1.1.3}}~~~~&=&
\sum_{i,j} \vartheta(c_{(1)}\ot a_{\psi\psi'}a_i) \lan a_i^*,a_{j\Psi}\ran
c_{(2)}^{\psi'\Psi}\ot a_j^*\ot c_{(3)}^{\psi}\\
&=&
\sum_{j} \vartheta(c_{(1)}\ot a_{\psi\psi'}a_{j\Psi}\ran
c_{(2)}^{\psi'\Psi}\ot a_j^*\ot c_{(3)}^{\psi}\\
{\rm \eqref{1.1.1}}~~~~&=&
\sum_{j} \vartheta(c_{(1)}\ot (a_{\psi}a_{j})_{\psi'}\ran
c_{(2)}^{\psi'}\ot a_j^*\ot c_{(3)}^{\psi}\\
{\rm \eqref{4.2.0}}~~~~&=&
\sum_{j} \vartheta(c_{(2)}\ot a_{\psi}a_{j}\ran
c_{(1)}\ot a_j^*\ot c_{(3)}^{\psi}\\
&=& c_{(1)}\ot \ol{\Omega}(c_{(2)}\ot a).
\end{eqnarray*}
Conversely, given $\ol{\Omega}$, we have to show that
$\alpha_1^{-1}(\ol{\Omega})=\vartheta$ satisfies \eqref{4.2.0}. Take any
$c\ot a\in C\ot A$ and write
$\ol{\Omega}(c\ot a)=\sum_l b_l^*\ot d_l\in A^*\ot C$. Since
$\ol{\Omega}$ is right and left $C$-colinear, we have
\begin{eqnarray*}
\ol{\Omega}(c_{(1)}\ot a_{\psi})\ot c_{(2)}^{\psi} &=&
\sum_l b_l^*\ot d_{l(1)}\ot d_{l(2)},\\
c_{(1)}\ot \ol{\Omega}(c_{(2)}\ot a)&=&\sum_l \lan b_l^*,a_{i\psi}\ran
d_{l(1)}^{\psi}\ot a_i^*\ot d_{l(2)}.
\end{eqnarray*}
Therefore we  can compute
\begin{eqnarray*}
\vartheta(c_{(2)}\ot a)c_{(1)}&=&
\lan \ol{\Omega}(c_{(2)}\ot a),1_A\ot\varepsilon\ran c_{(1)}\\
&=&\sum_l \lan b_l^*,a_{i\psi}\ran\lan a_i^*,1\ran \lan
\varepsilon,d_{l(2)}\ran
d_{l(1)}^{\psi}\\
&=& \sum_l \lan b_l^*1_{\psi}\ran d_l^{\psi}\\
&=& \sum_l \lan b_l^*1\ran d_l\\
&=& \lan \ol{\Omega}(c_{(1)}\ot a_{\psi}),1_A\ot\varepsilon_C\ran
c_{(2)}^{\psi}\\
&=& \vartheta(c_{(1)}\ot a_{\psi})c_{(2)}^{\psi}.
\end{eqnarray*}
Thus \eqref{4.2.0} follows. Finally, we show that $\alpha_1$ and
$\alpha_1^{-1}$ are  inverses of each other.
$$
\alpha_1^{-1}(\alpha_1(\vartheta))(c\ot a)=
\lan \vartheta(c_{(1)}\ot a_{\psi}a_i)a_i^*\ot c_{(2)}^{\psi},1_A\ot
\varepsilon_C\ran
= \vartheta(c_{(1)}\ot aa_i)\lan a_i^*,1_A\ran=\vartheta(c\ot a).
$$
We know that $\alpha_1(\alpha_1^{-1}(\ol{\Omega}))$ is right
$A$-linear. Hence suffices it to show that
$$\alpha_1(\alpha_1^{-1}(\ol{\Omega}))(c\ot 1)=c\ot 1,$$
for all $c\in C$. From \eqref{4.6.1}, we compute
$$\lan (a^*\ot c)b,1_A\ot \varepsilon_C\ran=
\lan a^*,b\ran \varepsilon(c)=\lan a^*\ot c,b\ot \varepsilon_C\ran.$$
Now write $\ol{\Omega}(c\ot 1)=\sum_r a_r^*\ot c_r$ and compute
\begin{eqnarray*}
\alpha_1(\alpha_1^{-1}(\ol{\Omega}))(c\ot 1)&=&
\lan \ol{\Omega}(c_{(1)}\ot a_i),1_A\ot \varepsilon_C\ran
a_i^*\ot c_{(2)}\\
&=&
\lan \ol{\Omega}(c_{(1)}\ot 1),a_i\ot \varepsilon_C\ran
a_i^*\ot c_{(2)}\\
&=& \ol{\Omega}(c\ot 1)_{[0]},a_i\ot \varepsilon_C\ran
a_i^*\ot \ol{\Omega}(c\ot 1)_{[1]}\\
&=& \sum_r \lan a_r^*,a_i\ran \lan \varepsilon,c_{r(1)}\ran a_i^*\ot
c_{r(2)}\\
&=& \sum_r a_r^*\ot c_r=\ol{\Omega}(c\ot 1).
\end{eqnarray*}
\end{proof}

\begin{theorem}\thlabel{4.9}
Let $(A,C,\psi)$ be a right-right entwining structure, and assume that
$A$ is finitely generated and projective as a $k$-module. With notation
as above, we have the following properties:\\
1) $F'$ is separable if and only if there exists $\ol{\Omega}\in V'_2$
such that for all $c\in C$,
$$\lan \ol{\Omega}(c\ot 1),1_A\ot\varepsilon_C\ran=\varepsilon_C(c).$$
2) $G'$ is separable if and only if there exists ${\Omega}\in W'_2$
such that for all $c\in C$,
$$\sum_i a_i(\varepsilon_C\ot I_A)\Omega(a_i^*\ot c)=\varepsilon_C(c)1.$$
3) The following assertions are equivalent:\\
a) $(F',G')$ is a Frobenius pair.\\
b) There exist $e\in W'_1$, $\vartheta\in V'_1$ such that $\Omega=
\beta_1(e)$ and $\ol{\Omega}=\alpha_1(\vartheta)$ are  inverses of each
other. \\
c) $A^*\ot C$ and $C\ot A$ are isomorphic objects in ${}^C{\cal M}(\psi)_A^C$.
\end{theorem}

\begin{proof}
We only prove $a)\Rightarrow b)$ in 3). First we show that $\Omega$ is
a left inverse of $\ol{\Omega}$. Since $\Omega\circ\ol{\Omega}$ is right
$A$-linear, it suffices to show that
\begin{eqnarray*}
\Omega(\ol{\Omega}(c\ot 1))&=&
\sum_i \vartheta(c_{(1)}\ot a_i)\Omega(a_i^*\ot c_{(2)})\\
&=& \sum_i \vartheta(c_{(1)}\ot a_i)\lan a_i^*,e^1(c_{(3)})_{\psi}\ran
c_{(2)}^{\psi}\ot e^2(c_{(3)})\\
&=& \vartheta(c_{(1)}\ot e^1(c_{(3)})_{\psi}\ran
c_{(2)}^{\psi}\ot e^2(c_{(3)})\\
{\rm \eqref{4.2.0}}~~~~&=&
\vartheta(c_{(2)}\ot e^1(c_{(3)})_{\psi}\ran
c_{(1)}^{\psi}\ot e^2(c_{(3)})\\
{\rm \eqref{4.4.3}}~~~~&=& c\ot 1.
\end{eqnarray*}
To show that $\Omega$ is a right inverse of $\ol{\Omega}$ we
use  that $\ol{\Omega}\circ \Omega$ is right $C$-colinear and
conclude that it suffices to show that for all $c\in C$ and $a^*\in A^*$,
$$(I_{A^*}\ot\varepsilon_C)(\ol{\Omega}(\Omega(a^*\ot c)))=
\varepsilon_C(c)a^*.$$
Both sides of the equation are in
$A^*$, so the proof is completed if we show that both sides
are equal when evaluated at an arbitrary $a\in A$. Observe that
$$\ol{\Omega}(\Omega(a^*\ot c))=\sum_i
\lan a^*,e^1(c_{(2)})_{\psi}\ran \vartheta\bigl((c_{(1)}^{\psi})_{(1)}\ot
e^2(c_{(2)})_{\Psi}a_i\bigr)a_i^*\ot (c_{(1)}^{\psi})_{(2)}^{\Psi}$$
hence
\begin{eqnarray*}
&&\hspace*{-15mm}
(I_{A^*}\ot\varepsilon_C)(\ol{\Omega}(\Omega(a^*\ot c)))(a)\\
&=& \lan a^*,e^1(c_{(2)})_{\psi}\ran \vartheta\bigl(c_{(1)}^{\psi}\ot
e^2(c_{(2)})a\bigr)\\
{\rm \eqref{4.3.2}}~~~~&=&
\lan a^*,(a_{\Psi}e^1(c_{(2)}^{\Psi}))_{\psi}\ran
\vartheta\bigl(c_{(1)}^{\psi}\ot e^2(c_{(2)}^{\Psi})\bigr)\\
{\rm \eqref{1.1.1}}~~~~&=&
\lan a^*,a_{\Psi\psi}e^1(c_{(2)}^{\Psi})_{\psi'}\ran
\vartheta\bigl(c_{(1)}^{\psi\psi'}\ot e^2(c_{(2)}^{\Psi})\bigr)\\
{\rm \eqref{1.1.3}}~~~~&=&
\lan a^*,a_{\psi}e^1((c^{\psi})_{(2)})_{\psi'}\ran
\vartheta\bigl((c^{\psi})_{(1)}^{\psi'}\ot e^2((c^{\psi})_{(2)})\bigr)\\
{\rm \eqref{4.4.4}}~~~~&=&
\lan a^*,a_{\psi}\ran \varepsilon(c^{\psi})=
\lan a^*,a\ran \varepsilon(c),
\end{eqnarray*}
as required.
\end{proof}

\section{The smash product}\selabel{5}
Let $(B,A,R)$ be a factorization structure (sometimes also called a
smash or twisted tensor product structure,
cf.\ \cite{Tambara}\cite[pp.\ 299-300]{Majid}\cite{Cap}).
This means that $A$ and $B$ are $k$-algebras
and that $R:\ A\ot B\to B\ot A$ is a $k$-linear map such that for all
$a,c\in A$, $b,d\in B$,
\begin{eqnarray}
R(ac\ot b)&=& b_{Rr}\ot a_rc_R,\eqlabel{5.0.1}\\
R(a\ot bd)&=& b_{R}d_{r}\ot a_{Rr},\eqlabel{5.0.2}\\
R(a\ot 1_B)&=& 1_B\ot a,\eqlabel{5.0.3}\\
R(1_A\ot b)&=& b\ot 1_A.\eqlabel{5.0.4}
\end{eqnarray}
We use the notation
$R(a\ot b)=b_R\ot a_R.$
$B\#_R A$ is the $k$-module $B\ot A$ with new multiplication
$$(b\# a)(d\# c)=bd_R\# a_Rc.$$
$B\#_R A$ is an associative algebra with unit $1_B\# 1_A$ if and only if
(\ref{eq:5.0.1}-\ref{eq:5.0.4}) hold. In this Section we want to examine
when $B\#_R A/A$ and $B\#_R A/B$ are separable or Frobenius. This will be
a direct application of the results in the second part of \seref{1a}.\\
In \seref{1a}, take $R=A$, $S=B\#_RA$. For $\ol{\nu}\in
V_1=\Hom_{R,R}(S,R)$, define $\kappa:\ B\to A$
by
$$\kappa(b)=\ol{\nu}(b\# 1).$$
Then $\ol{\nu}$ can be recovered form $\kappa$, since
$\ol{\nu}(b\# a)=\kappa(b)a$. Furthermore
$$a\kappa(b)=a\ol{\nu}(b\# 1)=\ol{\nu}(b_R\# a_R)=\kappa(b_R)a_R$$
and we find that
\begin{equation}\eqlabel{5.1.1}
V\cong V_1\cong V_3=\{\kappa:\ B\to A~|~a\kappa(b)=\kappa(b_R)a_R\}.
\end{equation}
Now we simplify the description of
$W\cong W_1\subset (B\#_RA)\ot_A (B\#_RA)$. Note that
there is a $k$-module isomorphism
$$\gamma:\ (B\#_RA)\ot_A (B\#_RA)\to B\ot B\ot A,$$
defined by
\begin{eqnarray*}
\gamma((b\# a)\ot (d\# c))&=& b\ot d_R\ot a_Rc,\\
\gamma^{-1}(b\ot d\ot c)&=& (b\# 1)\ot (d\# c).
\end{eqnarray*}
Let $W_3=\gamma(W_1)\subset B\ot B\ot A$. Take $e=b^1\ot b^2\ot a^2
\in B\ot B\ot A$ (summation implicitely understood). Then
$e\in W_3$ if and only if \eqref{1a.1.4} holds,
for all $s=b\# 1$ and $s=1\# a$ with $b\in B$ and $a\in A$,
if and only if
\begin{eqnarray}
bb^1\ot b^2\ot a^2&=& b^1\ot b^2b_R\ot a^2_R,\eqlabel{5.1.2}\\
(b^1)_R\ot (b_2)_r\ot a_{Rr}a^2&=& b^1\ot b^2\ot a^2a,\eqlabel{5.1.3}
\end{eqnarray}
for all $a\in A$, $b\in B$. This implies isomorphisms
\begin{equation}\eqlabel{5.1.4}
W\cong W_1\cong W_3=\{e=b^1\ot b^2\ot a^2\in B\ot B\ot A~|~
{\rm \eqref{5.1.2}~and~\eqref{5.1.3}~hold}\}.
\end{equation}
Using these descriptions of $V$ and $W$, we find immediately that
\thref{1a.2} takes the following form.

\begin{theorem}\thlabel{5.1}
Let $(B,A,R)$ be a factorization structure over a commutative ring $k$.\\
1) $B\#_RA/A$ is separable (i.e. the restriction of scalars functor
$G:\ {\cal M}_{B\#_RA}\to {\cal M}_A$ is
separable) if and only if there exists $e=b^1\ot b^2\ot a^2\in W_3$
such that
\begin{equation}\eqlabel{5.1.5}
b^1 b^2\ot a^2=1_B\ot 1_A\in B\ot A.
\end{equation}
2) $B\#_RA/A$ is split (i.e.
the induction functor $F:\ {\cal M}_A\to {\cal M}_{B\#_RA}$ is separable)
if and only if there exists $\kappa\in V_3$ such that
\begin{equation}\eqlabel{5.1.6}
\kappa(1_B)=1_A.
\end{equation}
3) $B\#_RA/A$ is  Frobenius (i.e. $(F,G)$ is Frobenius pair) if and only
if there exist $\kappa\in V_3$, $e\in W_3$ such that
\begin{equation}\eqlabel{5.1.7}
(b^2)^R\ot \kappa(b^1)_Ra^2=b^1\ot \kappa(b^2)a^2=1_B\ot 1_A.
\end{equation}
\end{theorem}

\thref{1a.7} can be reformulated in the same style. Notice that
$$\Hom_R(S,R)=\Hom_A(B\#_RA,A)\cong \Hom(B,A).$$
 $\Hom(B,A)$ has the following $(A,B\#_RA)$-bimodule structure
(cf. \eqref{1a.4.1}):
$$\bigl(cf(b\# a)\bigr)(d)=cf(db)a,$$
for all $a,c\in C$ and $b,d\in B$. From \prref{1a.5}, we deduce that
\begin{equation}\eqlabel{5.2.1}
V\cong V_2\cong V_4=\Hom_{A,B\#_RA}(B\#_RA, \Hom(B,A)).
\end{equation}
If $B$ is finitely generated and projective as a $k$-module, then we find
using \prref{1a.6}
\begin{equation}\eqlabel{5.2.2}
W\cong W_2\cong W_4=\Hom_{A,B\#_RA}(\Hom(B,A),B\#_RA).
\end{equation}
\thref{1a.7} now takes the following form:

\begin{theorem}\thlabel{5.2}
Let $(B,A,R)$ be a factorization structure over a commutative ring $k$,
and assume that $B$ is finitely generated and projective as a $k$-module.
Let $\{b_i,b^*_i~|~i=1,\cdots,m\}$ be a finite dual basis for $B$.\\
1) $B\#_RA/A$ is separable if and only if there exists an $(A,
B\#_RA)$-bimodule map
$\phi:\ \Hom(B,A)\cong B^*\ot A\to B\#_R A$ such that
$$\sum_i (b_i\# 1)\phi(b_i^*\# 1)=1_B\ot 1_A.$$
2) $B\#_RA/A$ is split
if and only if there exists an $(A, B\#_RA)$-bimodule map
$\ol{\phi}:\ B\#_R A\to \Hom(B,A)$ such that
$$\ol{\phi}(1_B\# 1_A)(1_B)=1_A.$$
3) $B\#_RA/A$ is  Frobenius if and only
if $B^*\ot A$ and $B\#_R A$ are isomorphic as $(A, B\#_RA)$-bimodules.
This is also equivalent to the existence of
$\kappa\in V_3$, $e=b^1\ot b^2\ot a^2\in W_3$ such that the maps
$$\phi:\ \Hom(B,A)\to B\#_R A,~~~\phi(f)=f(b^1)b^2\# a^2$$
and
$$\ol{\phi}:\ B\#_R A\to \Hom(B,A),~~~\ol{\phi}(b\# a)(d)=\kappa(bd_R)a_R$$
are  inverses of each other.
\end{theorem}

The same method can be applied to the extension $B\#_RA/B$. There are two
ways to proceed: as above, but applying the left-handed version of
\thref{1a.7} (left and right separable (resp. Frobenius) extension
coincide).
Another possibility is to use  ``op"-arguments. If $R:\ A\ot B\to B\ot
A$ makes
$(B,A,R)$ into a factorization structure, then
$$\tilde{R}:\ B^{\rm op}\ot A^{\rm op}\to A^{\rm op}\ot B^{\rm op}$$
makes $(A^{\rm op},B^{\rm op},\tilde{R})$ into a factorization structure.
It is not hard to see that there is an algebra isomorphism
$$(A^{\rm op}\#_{\tilde{R}}B^{\rm op})^{\rm op}\cong B\#_R A.$$
Using the left-right symmetry again, we find that
$B\#_R A/B$ is Frobenius if and only if
$(A^{\rm op}\#_{\tilde{R}}B^{\rm op})^{\rm op}/B$ is Frobenius
if and only if
$(A^{\rm op}\#_{\tilde{R}}B^{\rm op})/B^{\rm op}$ is Frobenius, and we can
apply Theorems \ref{th:5.1} and \ref{th:5.2}. We invite the reader to
write down explicit results.\\
Our final aim is to link the results in this Section to the ones in
\seref{2},
at least in the case of finitely generated, projective $B$.
Let $(A,C,\psi)$ be a right-right entwining structure, with $C$ finitely
generated and projective, and put $B=(C^*)^{\rm op}$. Let
$\{c_i,c_i^*~~i=1,\ldots,n\}$  be a dual basis for $C$. There is a
bijective correspondence between right-right entwining structures
$(A,C,\psi)$ and smash product structures $(C^{*\rm op},A,R)$. $R$ and
$\psi$ can be recovered from each other using the formulae
$$
R(a\ot c^*)= \sum_i \lan c^*,c_i^{\psi}\ran c_i^*\ot a_{\psi}, \qquad
\psi(c\ot a)= \sum_i \lan (c_i^*)_R,c\ran c_i\ot a_R.
$$
Moreover, there are
isomorphisms of categories
$${\cal M}(\psi)_A^C\cong {\cal M}_{B\# A}~~{\rm and}~~
{}_A{\cal M}(\psi)_A^C\cong {}_A{\cal M}_{B\# A}.$$ In particular,
$B\#_R A$ can be made into an object of ${}_A{\cal M}(\psi)_A^C$, and
this explains the structure on $C^*\ot A$ used in \seref{2}. Combining
Theorems \ref{th:2.10} and \ref{th:5.2}, we find that the forgetful
functor ${\cal M}(\psi)_A^C\to {\cal M}_A$ and its adjoint form a
Frobenius pair if and only if $C^*\ot A$ and $A\ot C$ are isomorphic as
$(A, (C^*)^{\rm op}\#_R A)$-bimodules, if and only if the extension
$(C^*)^{\rm op}\#_R A/A$ is Frobenius.

\end{document}